\documentclass[11pt]{article}
\usepackage{amsmath, amssymb}
\usepackage{amscd}
\usepackage{graphicx}
\usepackage{epstopdf}

\newtheorem{theor}{Theorem}[section]

\newtheorem{lem}[theor]{Lemma}
\newtheorem{propo}[theor]{Proposition}
\newtheorem{coro}[theor]{Corollary}

\newenvironment{pf}{{\it Proof.}}{\hfill $\square$\\}

\title{\bf   Harmonic functions on $\mathbb{R}$-covered foliations and group actions on the circle
\footnotetext{$2000$ \textit{Mathematics Subject Classification.} Primary 58J65, 37C85, 57R30; Secondary 53C12.} 
}

\author{S. Fenley\footnote{Partially supported by NSF grant DMS-0305313}, R. Feres, and K. Parwani}
\date{\today}
\begin{document}

\maketitle

\begin{abstract}
Let  $(M, \mathcal{F})$  be a compact  codimension-one foliated manifold whose leaves
are equipped with
 Riemannian metrics, and consider  continuous functions on $M$ that are harmonic along the leaves
   of $\mathcal{F}$. If every  such function is constant on leaves
we say that $(M,\mathcal{F})$ has the {\em Liouville property}.
  Our main result is that codimension-one foliated bundles over
  compact negatively curved manifolds satisfy the Liouville property. Related results
for $\mathbb{R}$-covered foliations, as well as
for discrete group actions and discrete harmonic functions, are also established.
\end{abstract}

\section{Introduction}

Let $M$ be a compact manifold 
and $\mathcal{F}$ a continuous foliation 
 of $M$ whose leaves are  $C^r$ Riemannian manifolds, $r\geq 2$. It is assumed 
 throughout the article
 that
 the boundary of $M$, if   non-empty, is a  union of  (compact) 
 leaves of
 $\mathcal{F}$. This implies that  all   compact leaves of $\mathcal{F}$  are closed manifolds.  
  The Riemannian 
 metrics on leaves, as well as their  derivatives up to order $r$, are assumed
  to vary continuously on $M$.
The pair $(M,\mathcal{F})$    refers here to foliations with
the given  choice of Riemannian metrics  even if the metrics are
not always  explicitly mentioned. The metrics yield
  Laplace-Beltrami operators on leaves varying continuously on $M$.

 Let $H(M,\mathcal{F})$ denote the set of
 real-valued functions on $M$ that are continuous 
 on $M$, $C^2$ on leaves,  and harmonic on leaves. We call 
 such functions 
   {\em leafwise harmonic}.
  If the leaves of $\mathcal{F}$ are Riemann surfaces, or more generally
  K\"ahler manifolds, we can similarly consider 
 the subset of $H(M,\mathcal{F})$ consisting of  the real part of   
{\em leafwise holomorphic} functions.   
   The continuous functions  that
  are constant on leaves, or {\em leafwise constant} functions,  form a subset of  $H(M,\mathcal{F})$. 
  If all leafwise harmonic (resp., holomorphic) functions are leafwise constant
  we say that
  $(M,\mathcal{F})$ 
has the  {\em Liouville} (resp., {\em holomorphic Liouville}) {\em  property}. 
The goal of this article is to study the Liouville property for  certain classes of
foliations.

 The problem of understanding which foliations
 have the Liouville property was first considered  in
 \cite{ghani1,ghani2}. A fairly detailed
 description of the structure of
 $H(M,\mathcal{F})$ in codimension $1$ under $C^1$ transversal
 regularity and in the absence of transverse invariant measures
 is obtained in \cite{abu}.  
 In order to provide some background for what will be proved here,
we briefly  list below a few pertinent results from these three papers.

 \begin{enumerate}
 \item For (real) codimension one foliations by K\"ahler manifolds (or more generally,
 foliations whose leaves are complex manifolds)  the holomorphic Liouville
 property holds. (\cite{ghani1}, Theorem 1.15.)
 \item  In \cite{abu} an example is given of 
 a codimension one foliation of a $3$-manifold  by Riemann surfaces
 for which the Liouville property does not hold.  The following  is also shown in \cite{abu} 
 (see Theorem 1.1 of \cite{abu} for the full statement): Let $\mathcal{F}$ be a $C^1$ 
 codimension one foliation of a compact manifold $M$ having no transverse invariant
 measures (in  particular, no compact leaves). Then there exists a finite number of minimal
 sets $\mathcal{M}_1, \dots, \mathcal{M}_k$ equipped with probability measures $\mu_1, \dots, \mu_k$
 such that   each 
 $f\in H(M,\mathcal{F})$ can be written uniquely  as a linear combination:
 $$f=\sum_{i=1}^k  \mu_i(f) \eta_i $$
  where the following notation is used: $\mu_i(f)=\int f d\mu_i$ and $\eta_i$ is a continuous, leafwise
  harmonic  function
  on $M$  which  gives
    the probability $\eta_i(p)$ that leafwise Brownian motion starting at $p$
  converges towards $\mathcal{M}_i$, for each $p\in M$.
  \item  It is shown in Theorem 4.1 of  \cite{ghani2} that there exists a foliated $S^2$-bundle over a compact Riemann surface, $(M,\mathcal{F})$,  such that:
  \begin{enumerate}
  \item The Liouville property does not hold for $(M,\mathcal{F})$;
  \item $(M,\mathcal{F})$ has exactly two minimal sets, $S_1$ and $S_2$, which are compact leaves homeomorphic to
  the base Riemann surface.  In the complement of $S_1\cup S_2$ the foliation and leafwise harmonic
  functions are smooth;
  \item The foliation is ergodic with respect to the smooth measure class. In particular, almost all leaves
  are dense in $M$.
  \end{enumerate} 
 \end{enumerate}
 
 The results of \cite{abu}, in particular item 2 above, 
are based on a study of the Lyapunov exponent for holonomy contraction
 along Brownian paths. They  depend in a crucial way on the foliation
 being $C^1$ and on the hypothesis that there are no transverse invariant measures.

 Given the above facts, particularly item 2,  it is natural to ask what can be said  about $H(M,\mathcal{F})$
   in
 codimension one when the results of \cite{abu} do not apply, namely when there are compact leaves present  or, more generally,
 transverse invariant measures, and/or the foliation is only $C^0$.
  In particular, we want to know under what natural
 hypothesis  codimension one foliations have the Liouville property.

  In this article we  restrict attention to $\mathbb{R}$-{\em covered}  or $I$-{\em covered} foliations. They are  
 defined by the property that the space of leaves of the induced foliation $(\widetilde{M}, \widetilde{F})$
   on the universal cover of $M$ is  Hausdorff. Equivalently, this leaf space is 
   homeomorphic to the real line or to the interval $I=[0,1]$,  respectively, hence the terminology.
	These are the simplest situations in terms of the topology of the foliation. In addition,
	as seen below, they exhibit a difficulty which is not covered by the results in \cite{abu}.
   Foliated circle bundles are $\mathbb{R}$-covered foliations. 
   Other examples of
   $\mathbb{R}$-covered foliations can
   be seen in \cite{fenley}.
 Based on what  we prove below it is natural to ask whether all 
 such  foliations have the Liouville property. As an initial 
support of an affirmative
 answer,
 we    mention   the following easy consequence of the topological structure
 of $\mathbb{R}$-covered foliations described in   \cite{fenley}.
 
 \begin{propo}\label{r}
 The Liouville property holds for 
   $\mathbb{R}$-covered  foliations
   without compact   leaves. 
    \end{propo}
 In fact, for foliations satisfying the conditions of proposition \ref{r}, we prove
 that every leafwise harmonic function is constant on $M$. 
 If there are no compact leaves, then we show there is only one minimal set, which
 then easily implies the Liouville property. Compare with results in \cite{abu},
 where one requires more than one minimal set to produce non trivial leafwise
 harmonic functions.

If there are compact leaves, the situation is much more interesting. First of all
 it is clearly  possible to have functions that are constant on leaves but not constant on
 $M$: when $\mathcal{F}$ is a fibration over the circle, any non-constant function
 on the circle pulls back to a leafwise constant, non-constant function on $M$. This 
 also happens to certain more general $\mathbb{R}$-covered foliations with compact leaves.
 
 Given proposition \ref{r}, our  problem is reduced to understanding what happens when there are compact leaves.
 In order to study leafwise harmonic functions or asymptotic behavior of holonomy (which is relevant
 here as well), it turns out that compact leaves are much trickier to understand. For example,
 the results of Deroin and Kleptsyn \cite{abu} do not apply when   there
 are compact leaves (even if one has the additional strong condition that holonomy is $C^1$).
 The same restriction 
 holds for the results of Thurston \cite{Th} on asymptotic behavior of holonomy.

   At this point we are not able to deal with the most general $\mathbb{R}$-covered foliations.
   For our main results we assume that
 the leaves of $\mathcal{F}$ have negative curvature $-$ this  is the condition 
 under which the 
 Liouville property might  be expected not to hold with greatest likelihood.
 Clearly, if the leaves of a foliation individually do not
 admit bounded, non-constant harmonic functions, then the foliated Liouville property  
 holds.  This is the case,
 for example, when the leaves are nilpotent covers of  recurrent (in particular, compact) Riemannian
 manifolds \cite{ls}, or  the Ricci curvature of leaves is non-negative \cite{schoen}.
 In dimension 3, results of Plante and Sullivan \cite{Pl,Su2} show that some form of negative curvature
 is the generic situation, at least in the large scale: if for example the leaves are 
  $\pi_1$-injective and $M$ is atoroidal and closed, then the leaves have negative curvature in
 the large, that is, they are Gromov hyperbolic. 
  In negative curvature, non-constant harmonic functions are plentiful,
  so if the Liouville property does hold  it must be due to
features pertaining to the transversal dynamics. 

Our main result is this:
\begin{theor}\label{main}
Let $(M,\mathcal{F})$ be a continuous codimension-$1$  foliated  bundle 
(with either circle or interval fibers)
over a
compact Riemannian manifold of negative sectional curvature. 
Then the Liouville property holds for $(M,\mathcal{F})$.  
\end{theor}
 
 By a {\em foliated bundle} $(M,\mathcal{F})$ we mean a foliation of the total space $M$ of a fiber
 bundle   whose fibers are everywhere transverse to the leaves of $\mathcal{F}$ and
 the local holonomy maps of the  fiber foliation are Riemannian
 isometries relative to the metric on the leaves of $\mathcal{F}$.

 Theorem \ref{main} is mainly a result about foliated interval bundles.  
The  claim for circle-bundles is an easy corollary given proposition
\ref{r}.  Just as easily,  theorem \ref{main}   implies the following: 
 
 \begin{theor}\label{moregeneral}
 Let $(M, \mathcal{F})$ be a continuous codimension-one foliation with negatively 
 curved leaves and let $\mathcal{M}$ denote (the closure of) the  union of all the minimal sets.
 Suppose that the metric completion
 $\widehat{U}$ of  each   component  $U$ of   $M\setminus \mathcal{M}$
 admits an interval-bundle structure   that makes the induced foliation
on  $\widehat{U}$ a   foliated interval-bundle over a compact base manifold.
Then the Liouville property
 holds for $(M,\mathcal{F})$.
 \end{theor}
 
 We note that a foliation satisfying the hypothesis of theorem \ref{moregeneral} is either
 minimal,  $\mathbb{R}$-covered, or $I$-covered foliation. 
 This can be seen as follows.
 A classical result of Haefliger states that the  union of all compact leaves
 of a codimension one  foliation is compact and there are finitely
 many compact leaves  up to isotopy in $M$. In addition, there are
 finitely many minimal sets which are not compact leaves.
  A proof 
 of this well-known fact
 for
 $C^2$ foliations can be found in \cite{candel}, theorem 8.3.2, and Cantwell
 and Conlon have a short, unpublished  proof for 
 $C^0$ foliations.
 So one possibility in theorem \ref{moregeneral} is that $\mathcal{F}$ is minimal,
 in which case the Liouville property holds by the maximum principle for harmonic 
 functions. 
 (A   continuous, leafwise harmonic function  must be constant on every minimal set. In fact, by
 the maximum principle the function
 is constant on
 a leaf where it attains its maximum or minimum value over a given  minimal set $\mathcal{A}$, hence it is constant on $\mathcal{A}$.) 
 Suppose now that  $\mathcal{F}$ is not minimal and let $\mathcal{A}$ be a minimal set.
 By the  just mentioned result of Haefliger's, if $\mathcal{A}$ is not a compact leaf and
 $B$ is a boundary leaf of $\mathcal{A}$, then $B$ is at a positive distance
 from any other minimal set.  Hence the hypothesis of
 theorem \ref{moregeneral} implies that if $U$ is a complementary component of $\mathcal{M}$ that has $B$ as one of its boundary leaves, then   $B$, and hence $\mathcal{A}$,  is a compact leaf.
 Therefore, we conclude that  the only minimal sets are compact leaves. 
 Given the
 finite number of isotopy  classes of compact leaves, it follows that by cutting
 $M$ along a compact leaf we obtain an $I$-bundle, and we can adjust the foliation to be
 transverse to the $I$-fibers. In particular,  the resulting  foliation is  $I$-covered.
 So  the original foliation (prior to cutting along a compact leaf)
  is either  $I$-covered or $\mathbb{R}$-covered.
 
 As an example to which 
 theorem \ref{moregeneral} applies, 
 start  with a foliated  interval-bundle and glue  the 
 boundary leaves with an arbitrary homeomorphism. This gives a foliation
 of a closed manifold satisfying  the conclusion of theorem  \ref{main}.
 Foliated circle bundles with compact leaves can be described in this way
 using a periodic  map as the gluing  map  so that  all points have the same
 period. To put things in perspective,  consider the situation in dimension $3$: the hypothesis
  of theorem \ref{main} implies that $M$ is Seifert fibered, the Seifert fibration
  given by circle fibers. In particular, it has a normal $Z$ subgroup. (See \cite{hempel}, chapter 12,
  for standard definitions.)  Theorem \ref{moregeneral},
  after cutting along a compact leaf, allows for any gluing  between top and bottom. The vast majority
   of such gluings yields hyperbolic $3$-manifolds. So this is much more general than theorem \ref{main}.

\vskip .15in
  Although   foliated bundles   may seem too restrictive a setting, they are   a very
  common type of foliation and the source of a large variety of examples and counter-examples
  in foliation theory. They are exactly the foliations that are associated with 
  group actions on the fiber space (the circle or interval, in theorem \ref{main}).
  The study of leafwise harmonic functions on codimension one foliated  bundles
  leads to interesting dynamical properties about group actions on $S^1$ or $I$. These are described now.
 
 Let $X_0$ denote the space of all harmonic functions $h$ on the unit open disc $\mathbb{D}$
 in $\mathbb{R}^2$ such that $|h(z)|\leq 1$ for all $z\in \mathbb{D}$.
 Then $X_0$ with the topology of uniform converge on
 compact subsets of $\mathbb{D}$
 is a compact metrizable space. Let  $\Gamma$ be a  group of hyperbolic isometries of the disc.
  $\Gamma$   
 acts on $X_0$ by composition: $ \gamma\cdot  h:= h\circ \gamma^{-1}$
 for $(\gamma, h)\in \Gamma\times X_0$. 
 The dynamics of this action can be    complicated even when 
 $\Gamma$ is only  an infinite cyclic group. For example, 
 it is shown in 
 \cite{ghani2}  that if $\Gamma$ is  cyclic   generated by 
  a parabolic or
 hyperbolic isometry of $\mathbb{D}$, the   
 action    admits
 a dense set of periodic orbits as well as orbits which are dense in $X_0$.

 It is  of   interest   to understand 
 what kinds of compact finite  dimensional  manifolds can arise as
 invariant subsets of $X_0$ for  a general $\Gamma$.
For example,  $S^2$   can,  but as we show below $S^1$ cannot.  
More precisely, there is an
 action of $\Gamma$ on $S^2$  
with respect to which one has a 
 non-trivial $\Gamma$-equivariant embedding  $S^2\rightarrow X_0$. 
This claim is essentially contained in \cite{ghani2}.
 Here, ``non-trivial'' means that
the   image of this map is not entirely contained  in the set of constant functions in $X_0$,
and   $F:X\rightarrow X_0$ from a given $\Gamma$-space $X$ into $X_0$
is said to be
$\Gamma$-equivariant if $F(\gamma(x))=\gamma \cdot F(x)$ for 
 all $\gamma \in \Gamma$ and $x\in X$.
The following is a corollary of  theorem \ref{main} when the base manifold
is a compact surface of constant negative curvature:

 \begin{coro}\label{circle1} Let $\Gamma$ be a discrete subgroup of  hyperbolic isometries
 of   $\mathbb{D}$ so that $\mathbb{D}/ \Gamma$ is a compact surface. 
 Consider an action of $\Gamma$ by homeomorphisms of  $X$, where $X$ is either $S^1$ or
  $[0,1]$.
 Then  
 any continuous, $\Gamma$-equivariant map  from $X$  into 
 $X_0$  takes values in
 the set of constant functions.
 \end{coro}

An action by homeomorphisms of the circle induces a foliated $S^1$-bundle
 over $\mathbb{D}/\Gamma$ by the suspension construction. A
  map from $S^1$ to $X_0$ as in the corollary produces a 
 function  on $S^1\times \mathbb{D}$ which is harmonic on   leaves and 
 induces a foliated harmonic function on the quotient $(S^1\times \mathbb{D})/\Gamma$
 by $\Gamma$-equivariance. 
 For the details of this easy proof see the general construction   in section
 5 of \cite{ghani2}. Similarly  for $I$ instead of $S^1$.
 Theorem \ref{main} then implies that 
 this function  is constant on   leaves, proving the corollary.

 We give now a somewhat different dynamical interpretation of the same result in the
 context of discrete harmonic functions.
Let    $\Gamma$ be, for the moment,  any countable group acting on a compact topological space
$X$ and equip $\Gamma$ with a probability measure $\mu$. 
Thus $\mu$
is a non-negative function on $\Gamma$ such that $\sum_{\gamma\in \Gamma} \mu(\gamma)=1$.
The  choice of $\mu$  specifies    transition probabilities
of  a random walk on $\Gamma$: the one-step transition 
from $\gamma$ to $\eta\gamma$   has probability  $\mu(\eta)$. 
To avoid trivialities, we assume that $\mu$ generates $\Gamma$;
i.e., the random walk starting from any $\gamma\in \Gamma$ has a positive
probability of reaching any other element of $\Gamma$ in a finite number 
of steps.
 We say that a continuous real-valued  function $f$ on $X$ is  {\em $\mu$-harmonic} if 
 $ f=P_\mu f$, where $P_\mu$ is the averaging operator defined by
  $P_\mu f(x)=\sum_{\gamma\in \Gamma} \mu(\gamma)f(\gamma(x))$
 for all $x\in X$.  The Liouville property in this discrete setting, for a given $\mu$, amounts 
 to all continuous, $\mu$-harmonic functions on $X$ being $\Gamma$-invariant.

Now suppose   that    $\Gamma$ is again a   group of isometries
 of $\mathbb{D}$ such  that $\mathbb{D}/\Gamma$ is a compact Riemann
 surface,  and let $\mu$ be a probability measure
 on $\Gamma$ that generates $\Gamma$.
  It makes sense to ask
 whether all actions    
 of $\Gamma$ on $X=S^1$ or $I$ by homeomorphisms have the Liouville property.
 This turns out to be true for at least one well-chosen $\mu$.
 In fact, 
 as first shown  by Furstenberg \cite{furs,ls,ancona}, there exists a probability measure $\mu$ on $\Gamma$
 with the following property: a bounded function on $\Gamma$ is $\mu$-harmonic
 for the action of $\Gamma$ on itself by left-translations 
 if and only if it is the pull-back to $\Gamma$   of a bounded harmonic
 function on $\mathbb{D}$ under the orbit map $\Gamma \rightarrow \Gamma\cdot z$, $z\in \mathbb{D}$.
 We call such  a measure a {\em discretization measure} on $\Gamma$.  In section  \ref{discretization}
 we give  a version of Furstenberg's  result
 for   the foliated bundle setting, theorem \ref{discrete}.
 Then theorem \ref{main}
 and theorem \ref{discrete}
 together imply the following corollary. The details are shown in 
 section \ref{discretization}.

 \begin{coro}\label{dyn}
 Let $\Gamma$ be a group of hyperbolic isometries of $\mathbb{D}$ such  that
  $\mathbb{D}/\Gamma$ is
  a compact surface  and let $\mu$ be
   a discretization   measure    on   $\Gamma$. For any
  representation  $\rho:\Gamma\rightarrow \text{Homeo}(X)$    of $\Gamma$ into the homeomorphism
 group of $X=S^1$ or $I$, 
 and   any continuous    $f:X\rightarrow \mathbb{R}$,  we have  
   $f\circ \rho(\gamma)=f \text{ for all } \gamma$ if and only if
   $P_\mu f =f$.
 \end{coro}
   
  Theorem \ref{discrete} also allows one to define a notion of discrete holomorphic
  function on a topological  $\Gamma$-space $X$, when $\Gamma$ is
  a cocompact group of isometries of a K\"ahler manifold.  A result 
 employing this idea is shown  in section \ref{holo}.

\vskip .15in
 We now give a brief sketch of the proof of theorem \ref{main}. 
 \begin{itemize}
 \item If   $\mathcal{F}$ has no
 compact leaves, then 
 $\mathcal{F}$ is a foliated circle bundle and it is
 $\mathbb{R}$-covered. Then the Liouville property is easily derived from the topological properties
 of such foliations discussed in section \ref{nocompactsect}, and properties of harmonic functions
 with respect to harmonic measures.
 \item If there are compact leaves,  we restrict attention to a connected component $U$ of $M\setminus \mathcal{K}$, the complement of the union of all compact leaves.  (Leafwise harmonic functions
 must be constant along   leaves in $\mathcal{K}$.)
 The metric completion of $U$
 is an interval bundle with compact boundary leaves and no compact leaf in the interior.
 This reduces the proof of the theorem to foliated interval bundles over a compact manifold and
 no interior compact leaves.
 These first 2 steps are done under the much more general condition of $\mathcal{F}$ being
 $\mathbb{R}$-covered or $I$-covered.
 \item The proof of the theorem for interval bundles proceeds by contradiction. We suppose
 that a nontrivial leafwise harmonic  continuous function $f$ exists, and normalize it so that
 it takes values $0$ and $1$ on the compact boundary leaves of $\mathcal{F}$. 
 Using the relationship between harmonic functions and properties of the foliated Brownian
 motion (under the assumption that leaves are negatively curved) we derive that $f$ is a monotone function on fibers of the interval bundle
 (lemma \ref{monotone}). 
 After blowing down interval bundles in $(M, \mathcal{F})$ where $f$ is constant
 along fibers, it  can be assumed that $f$ is  strictly monotone on fibers  (proposition \ref{minimal}).
 Both of these results make full use of the hypothesis of theorem \ref{main}: we need the foliated
 bundle property to directly relate  Brownian motion in different leaves. We also need negative curvature
 on the leaves to relate the harmonic function on the leaf with the behavior at infinity (this is
 done in the universal cover of the leaf).
 \item Using the strict  monotonicity  of  $f$ it is possible to define a new  foliated interval bundle
  topologically equivalent to the original one that is now {\em harmonic} in the following 
 sense: leaves of the new foliation are locally graphs of harmonic functions on the base manifold.
 This is shown at the beginning of section \ref{lipschitz}. Although the initial foliation was possibly
 only $C^0$, we prove in section \ref{lipschitz} that the new foliation is, in fact, Lipschitz continuous.
 \item In section \ref{harmonicmeasure} we prove the following general fact: If
 $M=K\times I$, where $K$ is a compact Riemannian manifold (no curvature assumption) and $I$ is the interval $[0,1]$,
  and $(M,\mathcal{F})$
 is a Lipschitz continuous harmonic foliation, then $\mathcal{F}$ is the product foliation. 
 This   result leads to a contradiction, since the original foliated interval bundle
 did not have compact leaves other than the boundary leaves.
\end{itemize} 

  The results of this article generate one obvious question: if $\mathcal{F}$ is $\mathbb{R}$-covered
  or $I$-covered,
  does $\mathcal{F}$ have the Liouville property? A key step to answering this question affirmatively
  is  to find
   some form
 of transversal monotonicity of leafwise harmonic  functions 
  as   obtained  in section 7 with the additional foliated $I$-bundle hypothesis. Another 
 very important question is 
 whether the curvature condition   can  be weakened.  In particular 
 what happens when the leaves are Gromov hyperbolic or negatively curved in the large, but
 not necessarily (Riemannian) negatively curved?

\section{Harmonic functions and  Brownian motion}\label{probability}
We begin by recalling some background material on harmonic functions and 
Brownian motion on Riemannian manifolds, with special attention to manifolds
of  negative sectional curvature. 
More details about Brownian motion 
 can be found, for example, in \cite{hsu}
or \cite{emery}. Brownian motion on foliated spaces is discussed in \cite{candelF} as well as
chapter 2 of \cite{candel2}. Information specific to negative curvature
can be found in \cite{ancona} and the other references to be cited below.

A few  key facts about harmonic functions   are listed below.
Let $N$ be a Riemannian manifold and $\Delta$ the Laplace-Beltrami operator on $N$.
A real-valued function $f\in C^2(U)$, where $U$ is an open set in $N$, is {\em harmonic}
on $U$ if $\Delta f=0$ on $U$. 
\begin{enumerate}
\item  The {\em maximum principle}:
 if $f$ is harmonic on a connected open set $U$ and attains a maximum 
 (or  minimum) value 
in $U$, then $f$ is  constant on $U$.
\item  The principle of {\em unique continuation} 
(see \cite{aronszajn} for a more general fact):
 if $f$ is harmonic on a connected
open  set $U$ and  constant on a neighborhood of some point in $U$, then $f$ is
constant on $U$.
\item The {\em Harnack inequality} \cite{moser}: If $U$ is open with compact closure
 and $V$ is a subset whose closure is contained in $U$, then there exists a constant $C>0$
 depending  only on $U$ and $V$ such that $\sup h|_V\leq C \inf h|_V$ for any  positive
 harmonic function $h$ on $U$.
 \item   The {\em Harnack principle} (see \cite{ancona},  p. 6): if $U$  is an open connected
 set in $N$ and $p\in U$, then the set of non-negative  harmonic functions $f$ on $U$ such that $f(p)=1$ 
 is compact in the topology of uniform convergence on compact subsets of $U$.     
\end{enumerate}

The standard probability setting for manifold-valued stochastic processes is assumed:
we fix throughout a probability space
  $(\Omega, \mathcal{B}, P)$ and
  a  {\em filtration} $\mathcal{B}_{*}=\{\mathcal{B}_t: t\geq 0 \}$
 of $\sigma$-algebras contained in $\mathcal{B}$. 
 That is,  $(\mathcal{B}_t)$ is 
 an increasing family of $\sigma$-algebras $\mathcal{B}_t$ containing  all sets of measure
 $0$ in $\mathcal{B}$. 
If $Y$ is an integrable real valued function on $(\Omega, \mathcal{B}, P)$, its expectation 
is denoted $E[Y]$, and if 
$\mathcal{A}$ a $\sigma$-algebra contained in $\mathcal{B}$, the 
  conditional
expectation of $Y$ given $\mathcal{A}$ is  denoted $E[Y|\mathcal{A}]$.
  Recall that a real valued stochastic process $\{Y_t : t\geq 0\}$ 
is a {\em martingale} if   $Y_t$  is integrable and $\mathcal{B}_t$-measurable for each  $t$
and for every pair $s, t\in [0, \infty)$, $s\leq t$, the equality   $Y_s=E[Y_t|\mathcal{B}_s]$
holds. A {\em Brownian motion} on  a Riemannian 
manifold $N$ with Laplace-Beltrami operator $\Delta$ is     an $N$-valued
stochastic process, $B_t$, $t\geq 0$, which  is continuous (i.e., sample paths $t\mapsto B_t(\omega)$
are continuous for a.e. $\omega\in \Omega$), adapted to the filtration (i.e.,
$B_t$ is $\mathcal{B}_t$-measurable for each $t\geq 0$), and for every smooth   function $f$
on $N$ the process
$$ M_t^f:= f(B_t) -f(B_0) -\frac12 \int_0^t \Delta f(B_t) dt$$
is a martingale. (This definition does not account for the possibility of explosions since
we will only deal with stochastically complete metrics later on.)
If 
$B_t$ is a Brownian motion on $N$ and if 
$\gamma: N\rightarrow N'$ is a local Riemannian isometry, then  $\gamma\circ B_t$
is   a Brownian motion on $N'$. For $p\in N$, Brownian motion conditional on 
$B_0=p$ will be written $B^p_t$. The corresponding conditional probability on $\Omega$
and expectation will be written
$P^p$ and 
 $E^p$, respectively.  
 Thus, for any bounded $f$ on $N$,  $E^p[f\circ B_t]:=
 \int_{\Omega} f(B_t(\omega))dP^p(\omega)$.

Let the Riemannian manifold  $N$ be geodesically  complete,  simply connected,
of  sectional curvature $K$  bounded by constants $-b^2\leq K \leq -a^2<0$.
 Let $S(\infty)$ be the sphere at infinity of $N$, which consists of equivalence classes
 of asymptotic geodesics. Then $\overline{N}=N\cup S(\infty)$
has a natural topology (the {\em cone topology}) that makes
 $\overline{N}$   compact and $S(\infty)$   its boundary. The latter is known as the {\em ideal boundary}
 of $N$.

 We   collect some of the main properties of
 Brownian motion on $N$ in the following list. (See  \cite{kifer} in addition to the references cited in
 each item.)

 \begin{enumerate}
 \item \label{b1} For any initial point $p\in N$,  $B_t^p$ converges in the cone topology, as $t\rightarrow \infty$,  
to   a  random point $B_\infty^p$  of  $S(\infty)$. (I.e., the path $B_t^p(\omega)$ converges
to a point $B_\infty^p(\omega)$ for $P^p$-a.e. $\omega\in \Omega$.)
\item \label{b2}
The probability   distribution of $B_\infty^p$  is a Borel probability measure $\mu_p$
 on $S(\infty)$.  Thus $\mu_p(A)$ is the probability  of the event
 $B^p_\infty \in A$, for a Borel $A\subset S(\infty)$. Its main property is that 
  $p\mapsto \mu_p(A)$ is a harmonic function on $N$.
 The measure $\mu_p$
  is known as the {\em harmonic measure}
 at $p$.
 This 
  should not be confused with harmonic measures for foliations as 
 defined by  Lucy Garnett in \cite{garnett}.
  The latter, which also plays a role in this paper, will be referred to
 either as   {\em stationary measures} for the foliated Brownian motion or  as   harmonic measures {\em in the sense of
 Garnett}.

 \item \label{b3}
 The measures $\mu_p$ are all equivalent among themselves.  This is a simple consequence
 of the maximum principle and that $p\mapsto \mu_p(A)$ is harmonic.  
 By the Harnack inequality, given 
 any pair of points $p, q\in N$ there exists  $C>0$   depending
 only on $p$ and $q$ such that 
 $ C^{-1}\mu_q(A)\leq  \mu_p(A)\leq C \mu_q(A)$
 for all $A$.
  The associated measure class defines
   the {\em harmonic measure class} of $S(\infty)$;

 \item \label{b4} For any bounded   function $g$ on $S(\infty)$, measurable relative to 
 the harmonic measure class,  the function
 $$H_g(p):=\int_{S(\infty)} g(\xi) d\mu_p(\xi)$$
 is harmonic on $N$. Conversely,
 if $H$ is a bounded  harmonic function on $N$,  there exists
 a bounded measurable  $g$ on $S(\infty)$, uniquely defined
 up to a set of harmonic  measure zero,  such that 
 $H=H_g$;

 
 \item \label{b5} If $H=H_g$ is a bounded harmonic function on $N$ with
 boundary values $g$, then $P^p$-almost surely $H(B^p_t)$ converges to
 $g(B_\infty^p)$ as $t\rightarrow \infty$;
 
 \item \label{b6} If $H=H_g$ is a bounded harmonic function on $N$ with boundary value
 $g$, then the non-tangential limit of $H$ exists almost everywhere on $S(\infty)$ with
 respect to the harmonic measure class. More precisely, 
 for $\xi\in S(\infty)$, $a>0$, and $t\mapsto r(t)$ a geodesic ray limiting at $\xi$,  
 denote  by
 $C_a(\xi)$ the set of all $p\in N$ such that the distance $d(p, r)<a$.
 Such a set is called  a {\em non-tangential cone} at $\xi$.
 Then, 
  for a.e. $\xi\in S(\infty)$ with respect to the harmonic measure class,
  and any non-tangential cone   $C_a(\xi)$, 
  $H(p)$ converges to $g(\xi)$ as $p\rightarrow \xi$   within $C_a(\xi)$. 
  (See \cite{anderson}.) 
 \end{enumerate}

 A note of caution: there is another natural probability measure 
on $S(\infty)$
 obtained by pushing forward the Lebesgue measure on the unit sphere
 $T^1_pN$ to $S(\infty)$ under the map that assigns  to each $v\in T_p^1N$
   the asymptotic class of the geodesic with initial condition $(p,v)$. 
  These measures  are known to be mutually equivalent for all $p$ and define
    the {\em geodesic measure class} on $S(\infty)$. 
  Even though the harmonic measures can be shown to be positive on 
  non-empty open sets  and to not have atoms \cite{kl},
     the geodesic  and the harmonic measure
   classes are in general mutually singular.
    In fact, by a result of A. Katok
 this is always the case for $N=\widetilde{K}$ and $K$ a closed surface
 of non-constant negative  curvature. If the sectional curvature is constant,
  the two measure classes coincide.
  
\section{Leafwise harmonic functions}   
  Let $M$ be a compact manifold and $\mathcal{F}$ a foliation of
$M$. Unless a stronger regularity assumption is explicitly
stated, $\mathcal{F}$ is a
continuous foliation with $C^2$ leaves. The tangent bundle
of $\mathcal{F}$ is given a Riemannian metric  smooth along leaves,
and the metric together with its derivatives of any order in the leaf direction
  are continuous 
in $M$.  We refer to this setting simply by saying that $\mathcal{F}$ is a foliation
of $M$ with Riemannian leaves.

The metric   induces a Laplacian on each
leaf of $\mathcal{F}$.
A continuous real valued function on $M$ is  {\em leafwise harmonic}
if  its restriction to each leaf   is
(smooth and) harmonic. 
Clearly,  a leafwise harmonic
function $f$ is constant on  any compact leaf, or on any leaf 
containing  a point of maximum or minimum value of $f$, due the the maximum principle. 
We say that a leafwise harmonic function is
{\em non-trivial} if it is not constant on at least one leaf of $\mathcal{F}$.

     Brownian motion on   leaves of $\mathcal{F}$ will still be denoted $B_t$.
Thus, for
 a  probability space  $(\Omega, \mathcal{B}, P)$ and 
 each $t\in [0,\infty)$, $B_t$ is a random variable  
 with values in $M$ and
   for $P$-a.e. $\omega\in \Omega$
 the path $t\mapsto B_t(\omega)$ is continuous and lies in the leaf of $\mathcal{F}$
 through $B_0(\omega)$.  
 The process and probability, conditional on beginning at $p\in M$, will be written as
 $B^p_t$ and $P^p$, respectively.

 \begin{propo}\label{constant1}
  Let $\mathcal{F}$ be a foliation of  a compact manifold $M$ with Riemannian leaves,
  as defined above.  Let 
  $L$ be a leaf of $\mathcal{F}$ with  sectional curvature
  $K_L$ satisfying at all points $-b^2 \leq K_L\leq -a^2<0$.  Let $S(\infty)$
  be the ideal boundary of the universal cover, $\widetilde{L}$,
   of $L$.
  Suppose that $f$ is a leafwise harmonic function on $(M, \mathcal{F})$ 
  and that the boundary values   of  the natural lift, $\tilde{f}$,  of  
  $f|_L$ to 
 $\widetilde{L}$
  are given by the Borel measurable function $g$ on $S(\infty)$.
  Then, for almost every $\xi\in S(\infty)$
  with respect to the harmonic class, there exists a leaf of $\mathcal{F}$ on
  which $f$ is constant and equal to $g(\xi)$.
      \end{propo}
\begin{pf}
Let $\xi\in  S(\infty)$ be a point of non-tangential convergence of $\tilde{f}$
and consider a sequence $p_n\in \widetilde{L}$ converging to $\xi$ along a geodesic ray.
For a fixed constant $c>0$, consider the sequence of balls $D(p_n,c)$ of
radius $c$ and center $p_n$. Then  for each $n$ and  all $q_n\in D(p_n,c)$ we have
 $\lim_{n\rightarrow \infty}\tilde{f}(q_n)=g(\xi)$.
 After passing to a subsequence,  the projection of $p_n$ to $M$ converges in $M$
 to a point $p$ and  the  balls converge to $D(p,c)$ as sets.
 Since $f$ is continuous on $M$, the value of $f$ on $D(p,c)$ is  
 equal to the  limit value $g(\xi)$. By the principle of unique continuation of harmonic functions
(see section \ref{probability}) $f$ must be constant, equal to $g(\xi)$, on that    leaf. 
\end{pf}

 \section{Foliations without compact leaves}\label{nocompactsect}
 
 If $\mathcal{F}$ is a foliation of a manifold $M$, let $\widetilde{\mathcal{F}}$
be the lift of $\mathcal{F}$ to the universal cover $\widetilde{M}$. The {\em space
of leaves} of $\widetilde{\mathcal{F}}$ is the quotient topological space   
$\widetilde{M}/\widetilde{\mathcal{F}}$ under the equivalence
relation that identifies points of $\widetilde{M}$ lying  on the same leaf.
A codimension one foliation $\mathcal{F}$ of a closed manifold $M$  is said to be  $\mathbb{R}$-{\em covered} (respectively, $I$-{\em covered}) if  
 the space of leaves, $\widetilde{M}/\widetilde{\mathcal{F}}$,  of  the  foliation 
 $(\widetilde{M}, \widetilde{\mathcal{F}})$  on  the universal cover of $M$  
is homeomorphic to $\mathbb{R}$ (respectively, to the closed interval $I=[0,1]$.)

\begin{propo}\label{topological}
Let $\mathcal{F}$ be an $\mathbb{R}$-covered foliation of a manifold $M$. 
Then one of
the three following cases happens:
\begin{enumerate}
\item There is a compact leaf;
\item $\mathcal{F}$ is minimal;
\item $\mathcal{F}$ is not minimal and there is a unique minimal set.
\end{enumerate}
\end{propo} 
\begin{pf} 
This is mostly contained in \cite{fenley}, proposition 2.6, although it is proved there
for
the special case of $3$-manifolds.

First suppose there are no compact leaves and let $Z$ be a minimal set. If $Z$ is all of $M$
we have alternative  2, so suppose this is not the case.
 We need to show that
$Z$ is unique. By lifting to a double cover we may assume that $\mathcal{F}$ is
transversely orientable. Let $U$ be a connected component of the
 complement of $Z$ and $\widehat{U}$ its metric completion.
Then $\widehat{U}$ has an {\em octopus decomposition} (proposition 5.2.14 of
\cite{candel}): $\widehat{U}=C\cup  A_1 \cup \dots \cup A_l$, where 
 $C$ is compact, $C\cap A_i$ is both  the transverse boundary of $A_i$
 and a connected component of the transverse boundary of $C$, and  the $A_i$ are
    $I$-bundles over non-compact
 manifolds $B_i$ and the foliation restricted to $A_i$ is transverse to the $I$-fibers. 
 The $B_i$ have boundary and the thickness  of the bundle goes
 to $0$ as distance from the   boundary of $B_i$ grows to infinity. 
 
 We claim that every leaf of $\mathcal{F}$ in $U$ has to 
 go  into some arm $A_i$ of $\widehat{U}$.
 In fact,
 let $D$ be a component of $\partial C \cap (\partial A_1 \cup \dots \cup \partial A_l)$.
 Then $D$ is contained in the transversal boundary of one of the $A_i$. Let $E$ be an $I$-fiber
  in $D$. Then $E$ connects $2$ horizontal boundary components of $A_i$.
 Lift $E$ to a transversal $\widetilde{E}$ in the universal cover connecting two boundary
 leaves of a connected lift $\widetilde{U}$ of $U$. Since the leaf space of $\widetilde{F}$
 is homeomorphic to $\mathbb{R}$, then the leaves in $\widetilde{U}$ all intersect $\widetilde{E}$.
 Projecting down to $M$ we obtain that all leaves in $U$ intersect $E$, hence $D$.
 
 The claim implies that every leaf of $\mathcal{F}$ in $U$ limits on points that the boundary leaves
of $A_i$ also limit on. 
This is because the thickness of the arms $A_i$ converges to zero as distance from the core
goes to infinity.
Therefore, any leaf in $U$ must limit on $Z$, hence it cannot be part of
another minimal set, proving the third alternative. 
\end{pf}

Notice that items 1 and 3 in proposition \ref{topological} are not mutually exclusive.

In Dippolito's work the $A_i$ are called foliated $I$-bundles. Here we restrict that terminology
to foliations by Riemannian leaves so that local holonomy along the $I$-fibers are Riemannian
isometries. (See the next section.)

Proposition \ref{topological} implies the Liouville property for
 $\mathbb{R}$-covered foliations  without compact leaves:

\begin{coro}\label{nocompact}
Let $(M,  \mathcal{F})$ be a compact foliated space with Riemannian leaves. Suppose that
the foliation is 
$\mathbb{R}$-covered without compact leaves. Then continuous,  leafwise harmonic functions are constant on $M$.
\end{coro}
\begin{pf}
Let $g$ be  continuous  leafwise harmonic.
By continuity, the closure of a leaf on which $g=c$, for some   constant
$c$,  contains a minimal set where $g=c$.
By the maximum principle, 
the maximum and minimum values of $g$ must be attained at points contained 
in leaves where $g$ is constant. If there are no compact leaves, the
previous proposition says that there is only one minimal set, therefore the maximum and minimum values  
of $g$ coincide. 
\end{pf}

There is more that can be said about $\mathbb{R}$-covered foliations of
a compact $M$ when there are compact leaves:
\begin{propo}\label{cut}
Let  the $\mathbb{R}$-covered foliation  $\mathcal{F}$ be transversely orientable and
have a compact leaf $K$.
Let  $T$ be the manifold obtained by cutting $M$ along $K$
 and letting $\mathcal{F}_1$ be the induced foliation on $T$.
Then $\pi_1(K)$ surjects in $\pi_1(T)$. If, in addition, $\dim M=3$
and $M$ is not doubly covered  by $S^2\times S^1$, then
$T$ is an $I$-bundle and $\mathcal{F}_1$ is isotopic to a foliation 
transverse to the $I$-fibers. 
\end{propo}
\begin{pf}
The claim about foliations in dimension $3$ can be found in the proof of lemma 2.5 of \cite{fenley}.
This uses the fact that the foliation in $T$ is $I$-covered and hence it is taut: every
two leaves are connected by a transversal to the foliation.
 The first claim can be seen as follows. 
 Let $\gamma$ be a loop in $T$ starting in $K$. Lift $K$, $\gamma$ to $\widetilde{K}$
 and $\tilde{\gamma}$ starting at $p$. By transverse orientability then $K$ locally separates
 $M$. Since the leaf space of $\widetilde{\mathcal{F}}$ is $\mathbb{R}$ it follows that
 $\widetilde{K}$ is the unique lift of $K$ to $\widetilde{T}$. Therefore $\tilde{\gamma}$
 ends in $\widetilde{K}$.  As $\widetilde{T}$  is simply
 connected, then $\tilde{\gamma}$ is homotopic to an arc in $\widetilde{K}$, so $\gamma$
 is homotopic to a loop in $K$. 
\end{pf}

Proposition \ref{cut} makes it clear, at least in dimension $3$, 
  that in trying to prove the Liouville property for
$\mathbb{R}$-covered foliations, it is essential to understand the case
of foliations transverse to $I$-fibrations. In the following sections we study 
the Riemannian  version of this, which we refer to as foliated $I$-bundles.

 \section{$I$-covered foliations}

We denote by $\mathcal{H}=\widetilde{M}/\widetilde{\mathcal{F}}$
the space of leaves  of the lifted foliation  to the universal cover $\widetilde{M}$ of
$M$.

\begin{propo}\label{limit}
Let $\mathcal{F}$ be a codimension-$1$ foliation of a compact manifold
 $M$   with boundary $\partial M= A_0\cup A_1$, where $A_0$ and $A_1$
 are leaves of $\mathcal{F}$. Suppose that 
 no leaf of $\mathcal{F}$ other than $A_0$ and $A_1$ is compact and
 that the space of leaves of   $\widetilde{\mathcal{F}}$
 is Hausdorff.
 Then the leaf space  $\mathcal{H}$ of $\widetilde{\mathcal{F}}$ is homeomorphic to a closed
 interval whose endpoints correspond to the unique lifts of ${A}_0$ and ${A}_1$,
 and   every  interior leaf   limits on both $A_0$ and $A_1$.
 \end{propo}
\begin{pf}
We first show that $\mathcal{H}$ is homeomorphic to $[0,1]$.
Suppose that there is a transversal arc in $\widetilde{M}$ connecting
a leaf of $\widetilde{\mathcal{F}}$ to itself. Join the endpoints by a path in the leaf to
produce a closed curve. Since $\widetilde{\mathcal{F}}$ is transversely orientable, this
path can be perturbed to produce a closed transversal, $\gamma$, to $\widetilde{\mathcal{F}}$.
As $\widetilde{M}$ is simply connected, $\gamma$ bounds a singular disc, $D$,
which can be assumed to be in general position with respect to $\widetilde{\mathcal{F}}$.
(See corollary  7.1.12 of \cite{candel}.)
In particular, $\widetilde{\mathcal{F}}$ is transverse to the boundary of $D$
and it induces on $D$ a singular $1$-dimensional  foliation, $\mathcal{F}^*$. The leaves 
of $\mathcal{F}^*$ are transverse to the boundary of $D$
and all singularities are isolated. By a standard argument
\marginpar{\tiny }
there must be a limit cycle, $\gamma$,  in $D$, and the germ of holonomy of $\mathcal{F}^*$
is contracting on at least one side of $\gamma$.
(See, for example, proposition 7.3.2 of   \cite{candel}; it is known that 
this argument, which is related to the Poincar\'e-Bendixson theorem,
can be carried out for $C^0$ foliations; see \cite{solodov,gabai}.)

 This closed curve lies on a leaf, $B$, of
 $\widetilde{\mathcal{F}}$ having contracting   holonomy germ along $\gamma$.
  But then, there are many leaves of $\widetilde{\mathcal{F}}$  near $B$
  which cannot be separated from $B$, contradicting the assumption that $\mathcal{H}$
  is Hausdorff.
  
  Hence any transversal to $\widetilde{\mathcal{F}}$  intersects a given leaf
  at most once, and so $\mathcal{H}$ is a $1$-manifold.
  It is clearly simply connected. In addition, it  has a countable basis  and is  Hausdorff
  by assumption. Therefore, $\mathcal{H}$ can only be homeomorphic to
  $(0,1)$, $[0,1)$, or $[0,1]$. But it has at least two boundary points, 
  which must be lifts of $A_0$ and $A_1$. In particular, $A_0$ and $A_1$ have
  unique lifts to $\widetilde{M}$, denoted
  ${A}'_0$
  and ${A}'_1$. It follows that $\mathcal{H}$ is homeomorphic to $[0,1]$,
  where $0$ and $1$ are identified with $A'_0$ and $A'_1$, respectively.  
  
  We now show that the interior leaves of $\mathcal{F}$ must limit on both
  $A_0$ and $A_1$. First observe that 
    $\mathcal{F}$ is transversely orientable. If not, some element of the fundamental 
    group of $M$ would switch the leaves  $A'_0$ and $A'_1$ in $\mathcal{H}$, and these would
    project to a single leaf in $M$, which is not the case.
    Suppose that an interior leaf $L$ does not limit on one of the boundary
    leaves, say $A_0$. Consider all the lifts of $L$ to $\widetilde{M}$.
    Each of them separates ${A}'_0$ from ${A}'_1$. This is because
    $\mathcal{H}$  is homeomorphic  to $[0,1]$ and the projection from $\widetilde{M}$
    to $\mathcal{H}$ is continuous, so a  path from ${A}'_0$ to ${A}'_1$ produces a path from 
    $0$ to $1$. 
    
    Let  $\mathcal{T}$ denote 
    the subset  of $\mathcal{H}$ corresponding to leaves of $\widetilde{\mathcal{F}}$ that 
    are separated from $A'_0$ by some lift of $L$. The above properties show
    that $\mathcal{T}$ is connected
  and homeomorphic to an interval $(c,1]$ or $[c,1]$.  Clearly 
  $c<1$ since any lift of $L$ separates $A'_0$ from $A'_1$ and also $c>0$, due
  to the assumption that   $L$ does not limit on $A_0$.  Let $\Theta$ be
  the projection  map from $\widetilde{M}$ to $\mathcal{H}$.
         Let $C$ be the   leaf of 
         $\widetilde{\mathcal{F}}$ corresponding to $c$. In particular $C$ is 
         not a lift of $A_0$ or $A_1$.
         We will show that $C$ projects to a compact leaf of $\mathcal{F}$, which
         is a contradiction.

 We claim that  any   covering translation of $\widetilde{M}$
 must map $C$ to itself. Covering transformations
 induce an action by orientation preserving homeomorphisms of
 $[0,1]$.
 If there is $h$ covering translation  so that $h$ does not leave $C$ invariant, then
 up to taking an inverse  we may assume that $h(c)<c$. 
  This contradicts the definition of $c$ as the infimum of $\Theta(V)$ where $V$
 is a lift of $L$.
 
  Let $\pi:\widetilde{M}\rightarrow M$ be the universal  cover projection. We now claim
  that $\pi(C)$ is compact. Otherwise there is a foliation box $Z$  in 
  $M$ in which 
   a sequence of distinct  sheets contained in
  $\pi(C)$ limit on a sheet of $\mathcal{F}$ in $Z$. 
  Lifting coherently to $\widetilde{M}$, we obtain
  a sheet $B'$ of $\widetilde{\mathcal{F}}$ and a sequence of distinct sheets in translates of $C$  that  converge to $B'$.
  But this was disallowed by the previous paragraph. This shows that $\pi(C)$ is compact,
   contradicting the hypothesis on $\mathcal{F}$.
\end{pf}

As an example to  which proposition \ref{limit} applies,
let $M=K\times I$, where $K$ is a compact Riemannian manifold and 
$I=[0,1]$, and  $\mathcal{F}$       a continuous foliation everywhere transverse
to the fibers of the fibration $\pi_2 : K\times I \rightarrow I$,    so that  
$A_i=K\times \{i\}$, $i=0,1$, are leaves of $\mathcal{F}$.
The proof of the previous lemma is much
simpler for this special case.

\section{Harmonic functions on $I$-covered foliations }
We describe here some  basic  consequences of assuming
 that 
  an $I$-covered foliation   carries 
  a non-trivial leafwise harmonic function.  
  Our goal is to show that under certain additional hypothesis there are no nontrivial
  such functions.

 \begin{lem}\label{function}
 Consider the same setting and assumptions of proposition \ref{limit}.
 If $(M, \mathcal{F})$ admits a non-trivial leafwise harmonic function, then
there exists a unique such function, $f$,  with the properties: the range of values
of  $f$ is the interval
 $[0,1]$; $f$  equals  $0$ on $A_0$ and $1$ on $A_1$; and
the restriction of $f$ to any leaf other than $A_0$ and $A_1$ has the range of
 values $(0,1)$. Any other   leafwise harmonic function $g$ is of the form $g=a f + b$
for constants $a, b$. 
\end{lem}
 \begin{pf}
 Let $g$ be a nontrivial leafwise harmonic function.
 As already remarked, $g$ is constant on each compact leaf, hence let 
   $a_0,  a_1$ be the constant values of $g$
 on   $A_0$ and $A_1$, respectively. Without loss of generality  we assume
 $a_0<a_1$. Note that  $a_1$ and $a_0$ are the maximum and minimum values of
 $g$ on $M$. In fact, suppose
 that  a maximum   value, $c$,  of $g$ was attained at an interior point, $q$. By the 
 maximum principle   the restriction of $g$ to the   leaf through $q$ 
  would be constant,  equal to $c$.
 Since
 an interior leaf must limit on both $A_0$ and $A_1$ by proposition
  \ref{limit}, then $a_0=a_1=c$. 
 This forces the maximum and minimum values of $g$ to coincide,  a contradiction. (Similarly,
 if $c$ is a minimum value.)
 The range
of $g$ on each interior leaf is the full open interval $(a_0, a_1)$
due, again, to   interior
leaves limiting on $A_0$ and $A_1$.
By composing $g$ with an affine function of the line we obtain $f$ with the claimed properties.
Uniqueness follows from the observation that if a leafwise harmonic function $h$ is $0$
  on $A_0$ and $A_1$, then the above argument shows that $h$ is identically
zero on all other leaves. 
 \end{pf}

For the next 2 results we assume that leaves of $\mathcal{F}$ have pinched negative curvature,
so we can use the facts of section \ref{probability}.

 \begin{lem}\label{01}
Assume, as in proposition \ref{limit},   that:
 $\mathcal{F}$ is a codimension-$1$ foliation of a compact manifold
 $M$   with boundary $\partial M= A_0\cup A_1$, where $A_0$ and $A_1$
 are leaves of $\mathcal{F}$;   
 no leaf of $\mathcal{F}$ other than $A_0$ and $A_1$ is compact; and
  the space of leaves of   $\widetilde{\mathcal{F}}$
 is Hausdorff. In addition, suppose 
that  the leaves of $\mathcal{F}$
  have negative sectional curvature
   and that $(M, \mathcal{F})$ admits a non-trivial leafwise harmonic function. 
  Let $f$ be the unique such function taking value $i$ on $A_i$,
  $i=0,1$.  Then the following properties hold:
 \begin{enumerate}
  \item $B_t$ is transient in $M\setminus \partial M$; that is, for any interior point $p$ of
 $M$, and any compact set $V\subset M\setminus \partial M$ containing $p$,  then
 for $P^p$-a.e. $\omega\in \Omega$, there is $\tau(\omega)<\infty$ such that  $B^p_t(\omega)$
 lies in the complement of $V$
for all $t\geq \tau(\omega)$. In other words, with probability one, $B_t$ converges towards 
$A_0$ or $A_1$;
 \item
  Let $L$ be the leaf through   $p\in M$, $\widetilde{L}$ the leaf
  through a lift  $p'$ of $p$, 
   $\tilde{f}$  the lift of $f$ to $\widetilde{L}$,
  and 
  $S(\infty)$ the ideal boundary of $\widetilde{L}$. 
    Then there exists a measurable set  $S_1\subset S(\infty)$
    for which the following holds:
 (i)  almost surely, Brownian motion
  $B^{p'}_t$  (in $\widetilde{L}$) converges to a point in $S_1$ if and only if   $B_t^p$ converges to $A_1$;  (ii) 
the probability that $B_t^p$ converges
  to   $A_1$ equals $f(p)=\mu_{p'}(S_1)$;   
  \item For every $p$ and a.e. unit vector $v\in T^1_p\mathcal{F}$
with respect to the harmonic class, viewed here as a measure class on
$T^1_p\mathcal{F}$, the geodesic ray with initial
  conditions $(p,v)$ converges to either $A_1$, if $v$ corresponds to
  $\xi$ in   $S_1$,  or  $A_0$ otherwise. (We make no similar claim for the
  visual measure on $T_p^1\mathcal{F}$.)
 \end{enumerate}
   \end{lem}
 \begin{pf}
 These assertions are 
   consequences of  proposition \ref{constant1}, lemma \ref{function},
  and the various facts about Brownian motion and boundary values of harmonic
  functions enumerated in section \ref{probability}. The curvature 
   pinching $-b^2\leq K\leq  -a^2<0$ holds since $M$ is compact.
  It is convenient to pass to the universal cover $(\widetilde{M}, \widetilde{\mathcal{F}})$.
 The lifts of the two compact leaves are denoted $A_i'$, $i=0,1$. 
  Then the above properties follow from the corresponding assertions for the lifted
  Brownian motion. 
  
  A key point to note is that, as $t\rightarrow \infty$, the distance between  $B_t^{p}(\omega)$  
and  $A_i$ goes to zero     if and only if $\tilde{f}(B_t^{p'}(\omega))$ converges to $i$ 
 since 
  ${f}$ is continuous,
 equals  $i$ on $A_i$, maps interior points of ${M}$ into $(0,1)$,
 and $\tilde{f}(B_t^{p'}(\omega))=f(B_t^p(\omega))$. This occurs because the
 corresponding fact holds in $M$, since $M$ is compact.

 Now the limit $\tilde{f}(B^{p'}_t)$ exists with $P^{p'}$-probability $1$ and equals $g(B_\infty^{p'})$,
 where $g$ is a function on $S(\infty)$ such that $\tilde{f}=H_g$. But by proposition
 \ref{constant1}, and since $\tilde{f}$ is not constant
 on any leaf except for $A_0'$ and $A_1'$,  it follows that, 
 almost surely, $g$   only takes the values $0$ and $1$.
Therefore, $g$ can be taken to be the indicator function of a subset of $S(\infty)$, denoted $S_1$.
 This shows assertions 1 and 2. The last statement of assertion 2 follows from property \ref{b4}
 of Brownian motion.
 Assertion 3 is a consequence of the existence of non-tangential limits of harmonic
 functions on $S(\infty)$. (Property \ref{b6} of section \ref{probability}.)
 \end{pf}

The  main conclusion  of lemma \ref{01} (parts 1 and 2) is summarized in the next corollary. 
 \begin{coro}\label{f}
Let
 $\mathcal{F}$ be a codimension one foliation of a compact manifold
 $M$   with boundary $\partial M= A_0\cup A_1$, where $A_0$ and $A_1$
 are leaves of $\mathcal{F}$;   
 no leaf of $\mathcal{F}$ other than $A_0$ and $A_1$ is compact;  
  the space of leaves of   $\widetilde{\mathcal{F}}$
 is Hausdorff; and leaves have negative sectional  curvature.
  If  $(M, \mathcal{F})$ admits a non-trivial leafwise harmonic function,
then the unique such function $f$ taking values $i$ on
 the boundary leaves $A_i$, $i=0, 1$, satisfies: $f(p)$ is the probability that the foliated 
 Brownian motion starting at $p$ converges to  $A_1$.
 \end{coro}
 
 We refer to the function $f$ as the {\em normalized} leafwise harmonic function on $(M,\mathcal{F})$.

 \section{Foliated bundles and monotonicity of $f$}
 It is natural to ask whether the normalized leafwise harmonic function $f$, which
 varies from $0$ to $1$ in the way from $A_0$ to $A_1$, is in some
 sense transversely monotone.
 It is not clear  how such a
property should be defined   for   general $I$-covered foliations, where
the manifold  may not even have an $I$-bundle structure.
Here we make the
additional restriction that $(M,\mathcal{F})$ be a foliated $I$-bundle, as defined below.

We first recall some definitions.
Let $K$ be a compact $n-1$-dimensional manifold  and $\pi :M\rightarrow K$ a fiber bundle
  whose fibers are everywhere transverse to a foliation  $\mathcal{F}$. 
  We assume that 
  the restriction 
  of $\pi$ to any leaf of $\mathcal{F}$ is a Riemannian covering of $K$.
  We say in this case that $(M, \mathcal{F})$ (together with the map $\pi$) is
  a {\em foliated bundle} with base manifold $K$.
  
  A foliated bundle also has the following description. Let $X=\pi^{-1}(q)$, $q\in K$,  represent a
  typical fiber of $\pi:M\rightarrow K$
  (a compact topological space) and let $\rho:\pi_1(K, q)\rightarrow \text{Hom}(X)$ denote the holonomy
  representation of 
  the fundamental group of $K$ acting on $X$ by homeomorphisms (or $C^r$ diffeomorphisms,
  if the foliation is $C^r$). Let  $\widetilde{K}$ be the universal covering of $K$. Then
  it can be shown that the quotient space $(\widetilde{K}\times X)/\Gamma$ for
  the natural action of  $\Gamma=\pi_1(K,q)$ on the product
  is isomorphic as a foliated bundle to  $(M, \mathcal{F})$.

  We are especially interested in the case where 
  $K$ has negative sectional curvature and   the fibers of the foliated bundle
  are homeomorphic to the interval $X=I=[0,1]$, where $0$ and $1$ are fixed
  points of $\rho$.
  We refer to  this setting
   as a {\em foliated $I$-bundle} with
  negatively curved leaves. For these $I$-bundles, $M$ has two boundary leaves, which are
  isometric to $K$. As already noted,  the foliation is transversely orientable
  since an orientation reversing  transformation  would have $0$ and $1$ in the same orbit
  of $\rho$, and
   $M$ would have only one boundary component rather than two.
  
  On $\widetilde{M}$, the map along $I$-fibers from $\widetilde{K}$ to the lift of
  any leaf is  a global isometry.
  Also $M$ is diffeomorphic to the product $K\times I$, so we can
  introduce a global {\em height  function} $\eta: M\rightarrow [0,1]$ corresponding to
  the projection on the second component of the product.
  This is a smooth function on $M$.
   Let $A_i$ be the boundary leaf of $M$
  corresponding to $\eta=i$, $i=0,1$.  
  Let $q$ be any point in $K$ and
  fix a lift $q'\in \widetilde{K}$. For any $p\in \pi^{-1}(q)$,
  let $L$ be the leaf of $\mathcal{F}$ through $p$. Then
   there is a unique
  local isometry $$\Phi_p:\widetilde{K} \rightarrow L,\text{ with } \Phi_p(q')=p \text{ and }
   \pi\circ \Phi_p:\widetilde{K}\rightarrow K$$ 
  is the universal
  covering map of $K$.

  Let $B^{q'}_t$ denote Brownian motion on 
  on $\widetilde{K}$ with  initial point $q'$. This is the same as the lift of 
  Brownian motion, $B^q$,  on $K$ with initial point $q$.
  Then   Brownian motion $B^p_t$ on $L$, for any $p$ in the
  fiber above $q$, 
  has a version given by $\Phi_p\circ B^{q'}_t$, which is also the lift
  to $L$ of $B_t^q$. This is because the restriction
   of $\pi$ to any leaf of $\mathcal{F}$
   is a Riemannian covering. The fact that Brownian motion along leaves
  can be, in this sense, ``synchronized'' along the $I$-fibers
  is 
  the main feature of the   Brownian motion
  on $(M,\mathcal{F})$ that we need here to deduce the
  property that if 
  a non-trivial leafwise harmonic function existed, then it would be monotone.
   This observation is the content of the next lemma.

\begin{lem}\label{monotone}
Suppose that the foliated $I$-bundle 
has no compact leaves other than $A_0$ and $A_1$, leaves have negative 
sectional curvature,   and
 there exists a non-trivial leafwise harmonic continuous
function.   For any $p\in M$, let
 $f(p)$ be  the probability that the foliated Brownian motion starting at $p$ will converge
 towards the boundary leaf $A_1$. Then for each $q\in K$, the restriction of
 $f$ to the fiber $\pi^{-1}(q)$ is a weakly monotone increasing function.
 \end{lem}
 \begin{pf}
 Recall that,
 if a non-trivial leafwise  harmonic function exists, 
  the
 unique such function taking values $i$ on $A_i$, $i=0,1,$ is $f$.
This is due to 
  corollary  \ref{f}.
 As $\mathcal{F}$ is transversely orientable, 
   given any two points $p_1, p_2$ in the fiber above $q$ 
   and any continuous curve
  $B^{q'}_t$ on $\widetilde{K}$ starting at $q'$, with $q'$ a lift of $q$ to $\widetilde{K}$, 
  we have
  $$\eta(p_1)<\eta(p_2)\Rightarrow  \eta(\Phi_{p_1}\circ B^{q'}_t)<\eta(\Phi_{p_2}\circ B^{q'}_t)$$
  for all $t\geq 0$.
  It follows that  the event $\Omega^1_{p_1}$ that $B_t^{p_1}$ limits on $A_1$
  can be regarded as a subset of the event $\Omega^1_{p_2}$ that $B_t^{p_2}$ limits
  on $A_1$. Since the probabilities of these events are $f(p_1)$ and $f(p_2)$, respectively,
  we must have $f(p_1)\leq f(p_2)$.
 \end{pf}

\begin{propo}\label{minimal} 
Let $(M,\mathcal{F})$ be a foliated $I$-bundle, $I=[0,1]$, with base manifold $K$,
where $K$ is a compact Riemannian manifold of negative sectional curvature.
We assume that there are no compact leaves in the interior of $M$ and that 
there exists a non-trivial leafwise harmonic function. Let $f$ be the normalized
such function.
 Then, after possibly
  blowing down interval-bundles in $(M, \mathcal{F})$, 
the restriction of $f$ to
each  $I$-fiber is a  strictly increasing function  onto $I$. 
\end{propo}
    \begin{pf}
     Let 
       $A_0$ and  $A_1$   be
as in lemma \ref{monotone}, and
     $U=M\setminus (A_0\cup A_1)$.
      We first make the following general observation. 
      Let $W$ be a noncompact foliated interval bundle in $U$. The lower and upper
      boundary leaves of $W$, denoted $L_0$ and $L_1$, respectively,
       are allowed to be the same. Let $p_i\in L_i$, $i=0,1$, be points in the same $I$-fiber.
      There is an isometry from $L_0$ to $L_1$ which sends Brownian motion in $L_0$ 
      starting at $p_0$ to Brownian motion  in $L_1$ starting at $p_1$.
           As in the proof of proposition \ref{topological}, consider the octopus decomposition of
      $W$.
   By lemma   \ref{01},
      Brownian motion $B_t^{p_0}$ in $L_0$ converges to either $A_0$ or $A_1$ almost surely. In particular, it
       escapes into
      the        arms of $W$ almost surely. 
   Since the thickness  of these arms  converges to $0$,
   then $B^{p_0}_t$ converges $A_1$ if and only if $B^{p_1}_t$ converges to $A_1$, and
   similarly for $A_0$. (Note that the index $i$ of the $A_i$ to which both 
   $B^{p_0}_t$ and $B^{p_1}_t$ converge is  a random variable, that
   is, a measurable function of the sample path.) Notice that
   this does not work in general if  $L_0$ and $L_1$ are not contained in a foliated
   $I$-bundle of $\mathcal{F}$.
 We now lift  all these sets  to the universal cover $\widetilde{M}$ of $M$ and
 let $\tilde{f}$ be the pull-back of $f$ to $\widetilde{M}$. Let 
  $\widetilde{L}_i$, $A'_i$,  denote the  lifts
   of $L_i$,  $A_i$, respectively, for $i=0,1$, where the $\widetilde{L}_i$ are
   boundary leaves of a connected lift of $W$.
  Consider the isometry $\Phi: \widetilde{L}_0
\rightarrow    \widetilde{L}_1$ defined via the holonomy map along   $I$-fibers and fix
   $p_i'\in \widetilde{L}_i$  such that  $\Phi(p'_0)=p_1'$.
Denote by   $S(\widetilde{L}_i)$ the set in the ideal boundary of $\widetilde{L}_i$
consisting of limit points, $B_\infty^{p_i'}$, of  Brownian paths converging
to $A'_1$. Then $S(\widetilde{L}_0)$
and $S(\widetilde{L}_1)$  are identified under the map induced by $\Phi$ on the ideal boundaries.
Therefore, $\tilde{f}|_{\widetilde{L}_0}$ and $(\tilde{f}|_{\widetilde{L}_1})\circ \Phi$
have almost surely the same boundary values at infinity
  and thus define the same harmonic function on $\widetilde{L}_0$.
   This  remark clearly  also applies to any pair of leaves between $\widetilde{L}_0$ and 
   $\widetilde{L}_1$.   This shows that $f$ is constant  along subsegments of $I$-fibers contained
   in $W$.

      First assume  that some leaf $L$ in $U$ accumulates only on $A_0$ and $A_1$.
      Let  $\widehat{W}$  be the metric completion of $W=U\setminus L$.  
      For every $p$ in $L$, the path starting at $p$ moving upwards along the 
      $I$-fiber of $p$ will hit $L$ again. But  $L$ does not limit in $U$ (that is, $L$ is
      properly embedded in $U$), so it makes sense
      to consider the first hit point from $p$ back in $L$. We obtain in this way a function
      from $L$ to itself
       that  is easily seen to be an isometry and is given by
      the holonomy map of an $I$-bundle structure on $\widehat{W}$. By the argument of the previous paragraph
      the restrictions of $f$ to the lifts of the boundary leaves of $\widehat{W}$ are equal
      on endpoints of $I$-fibers of $\widehat{W}$. But the top and
      bottom 
      boundaries of $\widehat{W}$ map  to the the same leaf $L$. Therefore,  the
      lifts of $W$ cover $\widetilde{M}\setminus (A'_0\cup A'_1)$. This shows that the
      restriction of $f$ to each $I$-fiber is constant, contradicting   that $f=i$ on $A_i$, $i=0,1$. 
      
      It will be assumed  from now on that every leaf in $U$ limits on $U$.  Suppose now that there 
      exist distinct points $p_0, p_1$ on the fiber $I_q$ of $q$ in the base manifold $K$ such that
      $f(p_0)=f(p_1)$. We want to show that these points lie  in the closure of a foliated
      $I$-bundle  in $U$. Take the interval $J\subset I_q$ with endpoints
      $p_0, p_1$ to be  maximal, i.e., $p_0$ is the lowest point in $I_q$ such that $f(p_0)=f(p_1)$,
      and  $p_1$ is the highest. Notice that $J$ is contained  in the
      interior  of $I_q$ because $f(0)=0$ and $f(1)=1$.
      Pass to the universal cover $\widetilde{M}$ and consider the
      harmonic function  $g=(\tilde{f}|_{E_1})\circ \Phi  - \tilde{f}|_{E_0}$, where  
       $E_i$ stands for  the leaf of $\widetilde{\mathcal{F}}$
      through lifts $p'_i$ of $p_i$ on the fiber $I_{q'}$ of a lift $q'$  of $q$, for $i=0,1$. Here, $\Phi$ is the fiber-respecting
      Riemannian isometry from $E_0$ to $E_1$ such that $\Phi(p_0)=p_1$.
      By  lemma \ref{monotone}, $g$ is a non-negative harmonic function on $E_0$ such that
      $g(p'_0)=0$.
      The maximum principle now implies that $g$ is identically $0$. 
      Let $J'$ be the interval  of $I_{q'}$ with endpoints $p_0', p_1'$.

      Now consider the   returns of $J$ to $I_q$ under the  foliation
       holonomy. From what has been shown, on any such interval return
       the function $f$ is constant and the interval is maximal relative to  this property.
       Therefore, the returns are either equal to $J$ or disjoint from $J$.
       In other words, the leaves of $\mathcal{F}$ through  $p_0, p_1$ are
       the boundaries of an $I$-bundle in $U$.  To see this, consider
       the set $W$ of leaves of $\widetilde{\mathcal{F}}$ through $J$.
       For any element $\gamma$ of $\pi_1(M)$, consider $\gamma(J)$. Move it
       along  by holonomy of $\widetilde{\mathcal{F}}$ to a subinterval $J_1$ of $I_{q}$.
       What has been shown above is that either $J_1$ equals $J$ or it is disjoint from $J$.
       This shows that $\gamma(W)$ is either equal to $W$ or disjoint from $W$. Hence
       $W$ projects to a foliated  $I$-bundle  in $M$ whose boundary leaves are
       the leaves through $p_0$ and $p_1$.

     So far we have  proved   that whenever $f$ takes the same
     value on two distinct points
    of  any $I$-fiber, 
     there is an $I$-bundle containing the two points such that $f$ is
     constant on each fiber of it, and the $I$-bundle is maximal
     relative to this property.

     Now let $\mathcal{U}$ be the union of the interiors of  all these $I$-bundles.
     The complement of $\mathcal{U}$ is a closed $\mathcal{F}$-saturated set in $M$
     and its intersection with $I_q$ is  a closed subset $V\subset I_q$    invariant
     under   holonomy of $\mathcal{F}$. $V$ does not have isolated points: if $v$ in $V$ is isolated
     then  the two open  intervals in $I_q\setminus V$ abutting $v$ would have 
     the same value of $f$, contradicting the maximality property above. Hence $V$ is a Cantor set.
     
      We can now collapse every $I$-bundle to one of its boundary leaves. This operation
      is done at most countably many times. There is an induced collapsed  foliation
       and induced continuous function   which is harmonic on leaves of the new foliation. The 
       new function  is now clearly strictly  monotone  along fibers.
        \end{pf}

 \section{The Lipschitz property}\label{lipschitz}
Although we have   assumed that the foliation and harmonic functions are only
continuous transversely, it turns out that more regularity can  be deduced in the case
of foliated $I$-bundles.
This fact will be essential in proving the main result.

 Let $(M,\mathcal{F})$ be as in  proposition \ref{minimal} 
 and suppose that a non-trivial, continuous leafwise harmonic function exists.
 Let $f$ be, as above, the unique   such function taking values $0$ and $1$ on the boundary components of $M$.  We assume that we have done the collapsing operation of proposition \ref{minimal} so that $f$ is strictly increasing on $I$-fibers.
 Let $\pi:M\rightarrow K$ be the bundle map and define
   $$\Psi:M\rightarrow K\times I,\  \ \Psi(p)=(\pi(p), f(p)).$$ Then $\Psi$ is a bijection. (Recall that if a fiber in $M$
 is identified with $I=[0,1]$ then the restriction of  $f$ to that fiber
  is a bijection from $I$ to itself.) Since $f$ is continuous on 
 $M$ and smooth along leaves, $\Psi$ is continuous on $M$ and smooth along leaves
 of $\mathcal{F}$. It can also be shown, using the strict monotonicity of $f$, that
 $\Psi$ maps $\mathcal{F}$ to a continuous foliation, $\mathcal{F}'$,  of $K\times I$ whose
 leaves are smooth and transverse to the fibers
 of the product fibration $\pi_1:K\times I\rightarrow K$. 
 We remark that $\mathcal{F}'$ is, like $\mathcal{F}$, a foliated bundle over the same base $K$, but
 it has the following additional  property: sheets of $\mathcal{F}'$ in any foliation box 
 of the form $\pi_1^{-1}(D)$
 are graphs of harmonic functions from  a sufficiently small 
 disc $D$ in $K$ to $[0,1]$. In particular,
if $S$ is an interval in a fiber $I_{q_0}=\pi_1^{-1}(q_0)$ and $S_q$ is the image
of $S$ under local holonomy of $\mathcal{F}$ from $I_{q_0}$ to $I_{q}$, then by fixing $q_0$
and letting $q$ vary,  the length of $S_q$ is a harmonic function of $q$,
since it is the difference between two locally defined harmonic functions corresponding 
to two different sheets of $\mathcal{F}'$.

Similarly,   the height function $\pi_2 :K\times I\rightarrow I$ is a non-trivial leafwise harmonic
 function for $\mathcal{F}'$.  A foliated $I$-bundle in  $K\times I$
  having the property that $\pi_2$ is
 leafwise harmonic will be called a {\em harmonic foliation}, and   
  $\mathcal{F}'$ may be viewed as a  ``harmonic straightening''
 of $(M, \mathcal{F})$. (No assumption on the curvature or topology of leaves is made
 in this definition.)
   It should be emphasized that the concept of a harmonic foliation
 is rather restrictive. In fact, under the fairly general assumptions of proposition
 \ref{product}, we show that a harmonic foliation is a product. 
 This will contradict  the existence of non-trivial
 leafwise harmonic functions under the conditions of  
 theorem \ref{main}.

  \begin{lem}\label{coordinatechange}
  Suppose that  the compact manifold $K$ has negative sectional  curvature
  and let $\mathcal{F}$ be a foliation of $M=K\times I$ having the properties:
  (i) $(M, \mathcal{F})$ is a harmonic foliation and (ii) no interior leaf is compact.
  Then $\mathcal{F}$ is transversely Lipschitz.
  \end{lem}
 \begin{pf}
 It is convenient to pass to the universal cover $\widetilde{M}=\widetilde{K}\times I$.
 The leaves of the lifted foliation, $\widetilde{\mathcal{F}}$, are now isometric to
$\widetilde{K}$ under the natural projection. We fix throughout the proof
two points  $q_1, q_2\in \widetilde{K}$ and use the natural parameter $t\in [0,1]$ to represent
 a point on the fiber $I_{q_i}=\{q_i\}\times [0,1]$,  $i=1,2$. 
 The holonomy map
 from $I_{q_1}$ to $I_{q_2}$ is then given by a strictly increasing function $t\mapsto H(t)$
 onto $[0,1]$. Our goal is to show that this function is Lipschitz.

 Due to  (i) and (ii), Brownian motion starting at  $p=(q, t)\in \widetilde{M}$
 converges to $A_1$ with probability $t$. 
 This follows from corollary \ref{f}
 and from the fact that  the height function, $h=\pi_2$, is  leafwise harmonic. (Clearly, this 
  also holds  for $t=0$ and $1$.)
 Let $L$ be the leaf of $\widetilde{\mathcal{F}}$
  through $p$ and  
  let $\mathcal{A}_1(L)$ be the measurable subset of the ideal boundary of $L$
  where $h|_L$ has boundary values equal to $1$. (Recall  lemma \ref{01}, part 2.)
  Then $$t=h(p)=\mu_p(\mathcal{A}_1(L)).$$ (This is due to  the same lemma; $\mu_p$ is
  the harmonic probability measure described in section \ref{probability}.)
  It will be convenient to be more explicit and denote 
  $\mathcal{A}_1(q,t):= \mathcal{A}_1(L).$ Note, however, that this set only depends on 
  the leaf $L$, so the following relation holds, by definition:
  $$\mathcal{A}_1(q_1,t)=\mathcal{A}_1(q_2,H(t)) $$
  for all $t\in [0,1]$. Using the Riemannian isometry $\pi_1|_L: L \rightarrow \widetilde{K}$
  we may identify the ideal boundary of any leaf $L$ with the ideal boundary, $S(\infty)$, of
  $\widetilde{K}$ and the harmonic measure $\mu_p$ with $\mu_q$, $q=\pi_1(p)$.
   From now on, we regard each $\mathcal{A}_1(q,t)$
  as a subset of $S(\infty)$ and write $t=\mu_q(\mathcal{A}_1(q,t))$.
  The same argument used to prove monotonicity of the leafwise harmonic
  function in lemma \ref{monotone} also shows that  
  $\mathcal{A}_1(q,s)\subset \mathcal{A}_1(q,t)$ whenever  $s<t$.
  This shows that, for any fixed $q\in \widetilde{K}$, 
  $$ t-s=\mu_q(\mathcal{A}_1(q,t)\setminus \mathcal{A}_1(q,s)).$$
  Now, for any given  Borel set $U\subset S(\infty)$, the  function  $q\mapsto \mu_q(U)$
  defined  on $\widetilde{K}$
   is  harmonic. By the Harnack inequality 
   there exists a constant $C=C(q_1, q_2)>0$ so that 
   $$ \mu_{q_2}(U)\leq C \mu_{q_1}(U)$$
 independent of $U$. It follows that
\begin{align*}
H(t)-H(s) &=\mu_{q_2}(\mathcal{A}_1(q_2,H(t))\setminus\mathcal{A}_1(q_2,H(s)))\\
                 &=\mu_{q_2}(\mathcal{A}_1(q_1,t)\setminus \mathcal{A}_1(q_1,s))\\
                 &\leq C \mu_{q_1}(\mathcal{A}_1(q_1,t)\setminus \mathcal{A}_1(q_1,s))\\
                 &=C(t-s).
 \end{align*}
 Therefore, $H$ is locally Lipschitz.  Now, 
 the Harnack  inequality  shows that
 the corresponding $C(q, q')$  is bounded
 for $(q, q')$ in a compact neighborhood of any $(q_1, q_2)$, hence $\mathcal{F}$
 can be covered by foliation charts with Lipschitz transition functions. 
 By compactness of  $M$  it follows that $\mathcal{F}$ is Lipschitz.
   \end{pf}

\section{Stationary measures} \label{garnett}
Stationary measures under foliated Brownian motion were introduced and studied in
\cite{garnett}, where they were named {\em harmonic measures}. 
Since we have been using the term to designate the measures $\mu_p$ on the
Poisson boundary (see section
\ref{probability}; this is the more traditional terminology from probability
theory)
we will refer to  Garnett's measures as {\em  stationary} (for the foliated  Brownian motion) or
 harmonic {\em in the sense of Garnett}.
See \cite{candelF}
and chapter 2 of \cite{candel2} for a comprehensive overview of the subject.

 The definition is as follows.
 Let  $(M, \mathcal{F})$ be a foliation by Riemannian leaves and
$\Delta$   the Laplace-Beltrami operator on leaves. 
A   Borel  measure $m$ on $M$
is harmonic in the sense of Garnett if $\Delta m=0$. By the duality between measures and
functions, this is interpreted by 
$\int_M \Delta \phi(x) dm(x)=0$, for all compactly supported smooth  functions, $\phi$, on $M$.
(By general measure theory, $m$ must be a regular Borel measure. See \cite{rudin}, theorem 2.18.)

It is shown in \cite{garnett} (see also proposition 2.4.10 of \cite{candel2}) that
  $m$   is harmonic in the sense of Garnett if and only  if,
on  any given foliated chart $U=D\times Z$ with transversal  $Z$,
$m$ can be disintegrated as $dm=h(q,t) d\sigma(q) d\nu(t)$, where
$\sigma$ is the measure on sheets induced by the Riemannian volume
form, $\nu$ is a measure on $Z$, and $q\mapsto h(q, t)$ is a non-negative harmonic function
on $D\times\{t\}$ for $\nu$-a.e. $t\in Z$.

\begin{propo}[Garnett]\label{constant}
Let $m$ be a harmonic probability measure on the foliated manifold
 $(M,\mathcal{F})$ and
$f$ a measurable, $m$-integrable, leafwise harmonic function on $M$. Then
$f$ is constant on $m$-a.e. leaf.
\end{propo}
\begin{pf}
We refer the reader to \cite{garnett} or \cite{candel2}. 
This corresponds to  proposition 2.5.6 of \cite{candel2}
and the fact that leafwise harmonic functions are precisely the
functions which are invariant under the diffusion semi-group, denoted $D_t$ in
\cite{candel2}.
\end{pf}

In the present paper, the function $f$ to which this proposition is applied is
continuous, so the conclusion is that $f$ is constant on any leaf in the support of
any harmonic measure in the sense of Garnett.  (Recall: 
a point  $p$ of a compact Hausdorff space lies in the support of a finite regular Borel
measure if, by definition, every neighborhood of $p$ has positive measure,
so arbitrarily close to $p$ in the support of $m$ there are leaves 
on which $f$ is constant.)

\section{Harmonic foliations and stationary measures} \label{harmonicmeasure}
In this section we use stationary measures  to prove that if $\mathcal{F}$
is a Lipschitz  harmonic foliation, then $\mathcal{F}$ is the product foliation.
This is then used in the next section to prove the main theorem.
The result of this section is very general in that we do not assume
that leaves of $\mathcal{F}$ have  negative curvature nor that $\mathcal{F}$
does not have compact leaves in the interior of $M$. These further conditions will be
imposed in the next section.

Recall the setup  of section \ref{lipschitz}:
Let $M=K\times I$, $I=[0,1]$, and let $\mathcal{F}$ be a    foliated $I$-bundle 
of $M$.  Let $h:=\pi_2:M\rightarrow I$ be the projection map. We assume that
$h$ is a leafwise harmonic function. In other words, the leaves of $\mathcal{F}$ are locally graphs of
harmonic functions. This  $(M,\mathcal{F})$ is called   
a harmonic foliation.  
In addition we assume  that $\mathcal{F}$ is Lipschitz continuous, that is,
its holonomy satisfies the Lipschitz property.

\begin{propo}\label{product}
Let $(M, \mathcal{F})$ be a Lipschitz continuous harmonic foliation of $M=K\times I$, where
$K$ is a compact Riemannian manifold. Then $\mathcal{F}$ is the product foliation.
\end{propo}
\begin{pf}
Let $\omega$ denote  the normalized Riemannian volume form on $K$, so that
the total volume is $1$, and let
$m$ denote the probability measure associated to the product volume form
 $\nu=\omega\wedge dt$  
on $M$.  
We claim that $m$ is a harmonic measure in the sense of Garnett.
As remarked in section \ref{garnett}, to show the claim  it suffices to verify  that 
the density functions   for the disintegration of $m$ on a foliation box
are harmonic on sheets.

  Let $\mathcal{W}$ be a foliation box of the form $\pi_1^{-1}(D)$,
where $D$ is a small enough Riemannian ball in $K$ with center $q_0$, and define the
the map $$\Phi:D\times I\rightarrow \mathcal{W} \ \text{ by }\ 
\Phi(q,t)=(q, \varphi(q,t)),$$ where for each fixed $t$, the graph of
$q\mapsto \varphi(q,t)$ is the sheet of $\mathcal{F}$ in $\mathcal{W}$ through
the point $(q_0, t)$. In particular, $\varphi(q_0,t)=t$. As $\mathcal{F}$ is a foliated
bundle,
the restriction of  $\pi_1$ to each leaf  is a local isometry onto $K$
 with local  inverse $q\mapsto (q, \varphi(q,t))$
for some $t$ and $q_0$. So  for a fixed $q_0$, we have that
$t\mapsto \varphi(\cdot, t)$ is a one-parameter family of isometries from $D$ to
sheets of $\mathcal{W}$.  Since the holonomy  of $\mathcal{F}$ is Lipschitz,
the function $\varphi$ is jointly  Lipschitz in $q$ and $t$ and smooth in $q$. 
Denoting by $I_q=\{q\}\times I$ the fiber $\pi_1^{-1}(q)$ above $q\in K$, note that
  the map $H_q(t)= \varphi(q, t)$ is the holonomy map over $D$ from $I_{q_0}$ to
  $I_q$. 
Thus $\Phi$ is a Lipschitz homeomorphism from $D\times I$ onto $\mathcal{W}$.

  As $\Phi$ is Lipschitz, 
  the pullback  $\nu'=\Phi^*\nu$
 is  a measurable, bounded form
  on $D\times I$ 
   by Rademacher's theorem on Lebesgue a.e. differentiability
  of Lipschitz functions. (Theorem 3.1.6 of \cite{federer}.)
   As $\nu'$ and $\nu$  are top-degree forms, 
  there is  a bounded measurable
  function $F$ on $D\times I$  such that
   $$\nu'_{q,t}=F(q,t) \omega_q \wedge dt.$$
   An elementary  Jacobian determinant calculation gives  that $F(q,t)=\varphi_t(q,t)$,
   whenever $\varphi_t(q,t)$ exists,   where
   $\varphi_t$ denotes partial differentiation with respect to $t$. In fact, for a.e. $(q,t)$ and
  all  $u\in T_qK\times \{0\}\subset T_{q,t}(K\times I)$,
   denoting $\tau=\frac{d}{dt}$, 
    then $d\Phi_{q,t}u=u+c \tau$
   for some scalar $c$  and $d\Phi_{q,t}\tau=\varphi_t(q,t)\tau$, so the
   determinant of $d\Phi_{q,t}$ with respect to a basis adapted to the product $K\times I$
    is   $\varphi_t(q,t)$.
   Furthermore,  the change of coordinates formula
   holds: 
  $$\int_{M}g\ \! dm= \int_0^1 \int_K (g\circ \Phi)(q,t) F(q,t) d\sigma(q)dt $$
   for any continuous function $g$ with compact support in $\mathcal{W}$.
   (This is easily derived from theorem 3.2.12 of \cite{federer}.)
   Here $\sigma$ is the normalized measure associated to the Riemannian
   volume form   $\omega$. In particular, this shows that the $F(q,t)$, when they exist,
   are the density functions for the disintegration formula in the foliation boxes.

    Therefore, to prove the claim that $m$ is harmonic we
    need to verify that $\varphi_t(q,t)$ is a harmonic function of $q$ on $D$ for
    almost every $t$. To see that this is the case we first define: 
  $$L(q,t, s)= \frac{\varphi(q,t+s)-\varphi(q,t)}{s}. $$
Note that $q\mapsto L(q,t,s)$ is a positive harmonic function for each fixed $t$ and $s$.
  Let $U$ denote the set of $(q,t)\in D\times I$ where 
  the limit of $L$ as $s\rightarrow 0$ exists. 
  This set is easily seen to be measurable.
  Since $\varphi$ is Lipschitz, $U$ has 
 full measure
  with respect to the product measure on $D\times I$.
  Now apply the standard Fubini theorem on product measure spaces
  to the characteristic function of $U$ to obtain that, since $U$ has full measure, 
   the slice 
  $$U_t:=U\cap (D\times \{t\})$$ has full measure for a.e. $t$. Therefore, 
  for a.e. $t\in I$,   $\lim_{s\rightarrow 0} L(q,t,s)$
exists (and is bounded) for a.e. $ q\in D$. 
If for a given $t$ the limit is $0$ for a.e. $q$, then $\varphi_t(\cdot, t)$
agrees a.e. with a (constant) harmonic function and the claim holds.
Now fix a $t$ for which $U_t$ has full measure in $D$ and
suppose that for some  $q'\in U_t$ the limit is positive.
Then the family $$l(\cdot ,t,s)=L(\cdot,t,s)/L(q',t,s)$$ of positive harmonic functions
satisfies $l(q',t,s)=1$ for each $s$.
By the Harnack principle (see section \ref{probability}), there is a subsequence
$s_n\rightarrow 0$ such that $l(\cdot, t, s_n)$ converges to a positive
harmonic function on $U_t$. Therefore $\varphi_t(\cdot ,t)$ agrees a.e. in $D$
with a positive harmonic function. This concludes the proof that $m$ is a harmonic
measure  in the sense of Garnett.

We can now apply proposition \ref{constant} to obtain that the height function
$h$ must be constant a.e. on the support of $m$. But $m$ has full support, 
so $h$ is constant on leaves everywhere. Therefore, $\mathcal{F}$ is the trivial foliation.
\end{pf}

\section{End of proof of the main theorem}
  
  We have finally obtained the desired 
   contradiction to the existence of non-trivial harmonic functions on foliated $I$-bundles.

  \begin{lem}\label{Ibundle}
 Let $(M,\mathcal{F})$ be a foliated $I$-bundle with base manifold $K$, where $K$ is
 a compact Riemannian manifold of negative sectional curvature. Suppose that
 no interior leaf is compact. Then $\mathcal{F}$ does not
 admit  continuous,  non-trivial,  
 leafwise harmonic functions. 
  \end{lem}
  \begin{pf}
  This is now an immediate consequence of propositions \ref{minimal},
  \ref{product} and lemma \ref{coordinatechange}. 
  Assume by way of contradiction that $(M,\mathcal{F})$ has a nontrivial
leafwise harmonic function $f$, which we take to be normalized.
Proposition  \ref{minimal} shows that after a blow down of $I$-bundles of $\mathcal{F}$ 
we can assume that $f$ is strictly increasing along $I$-fibers. The map $\Psi(p)$ defined
prior to  lemma \ref{coordinatechange} transforms this into a harmonic foliation in $K\times I$.
 Lemma \ref{coordinatechange} shows that this foliation is Lipschitz. Proposition \ref{product}
 shows that the new foliation  is a product foliation, contradicting the fact 
 that the original foliation did not have compact leaves in the interior.
    \end{pf}

 \begin{theor}  
 Let $(M, \mathcal{F})$ be a continuous foliated $S^1$-bundle with base manifold $K$, 
 where $K$ is a compact Riemannian manifold of negative sectional curvature. 
 Then $\mathcal{F}$ does not admit
 a non-trivial,  continuous, leafwise harmonic function.
 \end{theor}  
 \begin{pf}
 If there are no compact leaves, the theorem  reduces to corollary \ref{nocompact},
 since a foliated $S^1$-bundle is  $\mathbb{R}$-covered.
 Otherwise, 
 let $\mathcal{K}$ be the union of all compact leaves.
By a well-known theorem of Haefliger $\mathcal{K}$ is a compact set.
The maximum principle implies that
 any leafwise harmonic function on $M$ is leafwise constant on $\mathcal{K}$,
 so   we may assume that  $M\setminus \mathcal{K}$ is non-empty.
  Let $U$ be
 a component of the complement of $\mathcal{K}$.
 Then the metric completion of $U$ is an interval bundle with 
 compact boundary leaves
 and no compact leaf in the interior.
 We can now apply lemma \ref{Ibundle} to conclude the proof.
 \end{pf}
  
  Essentially the  same argument shows the  more general
  theorem \ref{moregeneral}. If there are no compact leaves, the result follows from
  corollary  \ref{nocompact}. Otherwise, the proof reduces to the foliated interval bundle
  case just as was done above for foliated circle bundles.
  
  \section{Discretization}\label{discretization}
  Let $\Gamma$ be a countable group of isometries of 
  a connected  Riemannian manifold $D$  such
  that $D/\Gamma$ is a compact manifold.
  We assume that $D$ is {\em transient}, i.e., for any $p\in D$, Brownian motion starting
  at $p$ eventually escapes any compact set almost surely. Of particular interest for
  us is  the hyperbolic disc $D=\mathbb{D}$. In addition to the properly discontinuous
  action on $D$ by isometries,   we assume that $\Gamma$ acts via homeomorphisms
  on a compact space $X$. For simplicity, the same notation will be used for both actions.

 We describe in this section a bijective correspondence between leafwise harmonic
 functions on the  foliated $X$-bundle $(M,\mathcal{F})$ over 
  $D/\Gamma$ and harmonic functions in a discrete sense
to be defined below for the   $\Gamma$-action on $X$.   
The reason for assuming   $D$  transient  is that, in the alternative
($D$  recurrent), bounded harmonic functions on $D$ are constant  
and the results below become trivial. (See, e.g., theorem 2.1, section 4,  of \cite{ancona}.)

  Let $V$ be a countable set and $P:V\times V\rightarrow [0,1]$ a
  Markov transition kernel. This means that 
  $\sum_{v\in V}P(u,v)=1$ for each $u\in V$.
  We regard $V$ as the set of states of a Markov chain with  probability
    $P(u,v)$ of transition from state $u$ to state $v$. 
    A real valued function $\varphi$ on $V$ is called $P$-{\em harmonic}
    if $\varphi=\mathcal{P}\varphi$, where we define
    $$\mathcal{P}\varphi(u)=\sum_{v\in V}P(u,v)\varphi(v)$$
    for each $u$. The transition probabilities can also be expressed by 
    a family of probability measures, $u\mapsto \mu_u$ on $V$, where
    $\mu_u(A)=\sum_{v\in A}P(u,v)$.

 Now  take $V$ to  be the orbit $\Gamma\cdot p_0$ of a point $p_0$ in      $D$.
  If 
  $P(\gamma u, \gamma v)=P(u,v)$ 
  for all $u,v\in V$ and $\gamma\in \Gamma$
    we say   that the Markov kernel 
  is {\em compatible} with the action of $\Gamma$ on $D$. More generally, it will be considered below
  functions $P:D\times V\rightarrow [0,1]$ such that $\sum_{v\in V} P(p, v)=1$ for all $p\in D$.
  We also refer to such $P$ as a Markov kernel and define compatibility similarly.
 
 Let $H_b(V, P)$ denote
 the space of all bounded  $P$-harmonic   functions on $V$
 and $H_b(D, \Delta)$ 
  the space of all bounded harmonic  functions on $D$ with respect to the
 Laplace-Beltrami operator.
 The following theorem
 says, in essence, that these spaces are isomorphic
 when $P$ is a  well-chosen
 Markov kernel   on $V$. The isomorphism amounts to restricting functions
 of $D$ to $V$.
In particular, 
   bounded harmonic functions on $D$ 
 can be completely recovered  given their values on only a discrete set of points in $D$.
  Theorem \ref{fls} is
   a special case, sufficient for our needs, of a discretization 
 property first observed by Furstenberg \cite{furs} for the group of isometries of $\mathbb{D}$
 and later generalized by Lyons and Sullivan in \cite{ls} and others.  We refer the
 reader to theorem 1.1,  section 4 of 
 \cite{ancona}.

 \begin{theor}[Furstenberg, Lyons-Sullivan]\label{fls}
 Let $D$ and $V=\Gamma\cdot p_0$ be as above. 
 There exists a Markov transition kernel $P:D\times V\rightarrow \mathbb (0,1)$ (strictly positive)
  such that
  the operation $\mathcal{P}:H_b(V,P)\rightarrow H_b(D, \Delta)$ defined
 by $$\mathcal{P}\varphi(p):=\sum_{v\in V}P(p,v)\varphi(v)$$ is a bijection.
 The inverse $\mathcal{P}^{-1}$ is the restriction operation $f\mapsto f|_V$. 
 Furthermore, $P$ is 
 compatible with the $\Gamma$-action on $D$ and $V$. 
 \end{theor}

The Markov kernel $P$ of theorem \ref{fls} is associated to a probability measure
$\mu$ on $\Gamma$ defined by $\mu(\gamma)=P(\gamma p_0, p_0)$. We call
$\mu$
a {\em discretization measure} on 
  $\Gamma$.  This is the measure referred to in corollary \ref{dyn}. 
Compatibility of $P$ with the $\Gamma$-action  implies $P(\gamma p_0, \eta p_0)=\mu(\eta^{-1}\gamma)$ for all $\eta, \gamma\in \Gamma$.

Our goal is to derive for foliated bundles a  discretization result similar to theorem \ref{fls}.
Let $X$ be a compact topological space and suppose that the group $\Gamma$ 
 of
isometries of $D$  also
acts on $X$. The latter action is given by an arbitrary   homomorphism of
$\Gamma$ into the group of homeomorphisms of $X$.

A probability measure $\mu$ on $\Gamma$ (shortly it will be assumed that
$\mu$ is a discretization measure) induces a Markov transition kernel
on $X$ by setting $P(x, y)$ equal
to the sum of    $\{\mu(\gamma) :\gamma\in \Gamma \text{ and } 
 y=\gamma x\}$, and $P(x,y)=0$ if $y$ and $x$ do not lie on the same $\Gamma$-orbit.
 Notice that for any given $x$ the probability $P(x,y)$ is nonzero for at most countably many $y$.
 (By a Markov transition kernel on $X$ we mean that the sum of $P(x,y)$ over   $y$  equals $1$
 for all $x$.)
The associated  operator $\mathcal{P}$ acting 
  on  continuous functions on $X$ is
$$ \mathcal{P}\varphi(x)=\sum_{\gamma\in \Gamma} \mu(\gamma)\varphi(\gamma x).$$
A function $\varphi$ on $X$ is said to be $\mu$-harmonic if $\varphi=\mathcal{P}\varphi$.
The space of  continuous $\mu$-harmonic functions on $X$ will be denoted by 
$H(X,\Gamma, \mu)$.

Let now $(M,\mathcal{F})$ be the   foliated $X$-bundle
associated to the given $\Gamma$-action on $X$. Thus $M=(D\times X)/\Gamma$ is the orbit space
for the action $\gamma(p,x)=(\gamma p, \gamma x)$. 
As before, the space of continuous leafwise harmonic functions on $M$ will be written
$H(M,\mathcal{F})$. We wish to define a sort of restriction map, $R$,
from the space of continuous functions on $M$ into the space of continuous functions on $X$.
Given  $f:M\rightarrow \mathbb{R}$ continuous,  write
$\tilde{f}:=f\circ \pi$, where $\pi : D\times X\rightarrow M$ is the projection map.
Notice that
$\tilde{f}(\gamma p, \gamma x)=\tilde{f}(p, x)$. 
Now fix a point $p_0\in D$ and define   $R$ by
$$Rf(x):=\tilde{f}(p_0, x).$$ 
 Let  $V=\Gamma \cdot p_0$, a discrete subset
of $D$.
We emphasize that the discretization measure $\mu$ in the next  theorem is the same one obtained from 
theorem \ref{fls} and does not depend on the choice of $\Gamma$-space $X$.

\begin{theor}[Discretization] \label{discrete} 
Let $D$ be a transient Riemannian manifold, $\Gamma$ a group of isometries of $D$
such that $D/\Gamma$ is a compact manifold, and
$\mu$  a discretization probability measure on $\Gamma$.
Let $(M,\mathcal{F})$ be a foliated bundle with fiber $X$and base space $D/\Gamma$.
Then the  restriction map   $R:H(M, \mathcal{F})\rightarrow H(X, \Gamma, \mu)$
is a bijection.
Furthermore, continuous  leafwise constant functions on $M$ correspond bijectively
under $R$ to $\Gamma$-invariant continuous 
functions on $X$.  In particular, $(M,\mathcal{F})$ has the Liouville property if and
only if continuous $\mu$-harmonic functions on  $X$ are $\Gamma$-invariant. 
\end{theor}
\begin{pf}
Having fixed $p_0\in D$, 
the measure $\mu$ is defined by $\mu(\gamma):=P(\gamma p_0,p_0 )$,
where $P$ is a $\Gamma$-compatible Markov kernel on $V=\Gamma\cdot p_0$
given by theorem \ref{fls}.  
 This implies that the condition $\hat{f}(u)=\sum_{v\in V} P(u,v)\hat{f}(v)$
  characterizing  a $P$-harmonic function $\hat{f}$  on $V$ is equivalent to
 \begin{equation}\label{a} \hat{f}(\gamma p_0)=\sum_{\eta\in \Gamma} \hat{f}(\eta p_0)\mu(\eta^{-1}\gamma)\end{equation}
 for all $\gamma, \eta \in \Gamma$. This uses the fact that $P(\gamma p_0, \eta p_0)=\mu(\eta^{-1}\gamma)$.
 
 Now, let $f\in H(M,\mathcal{F})$ and define the notation  $\Phi_{p_0}(x):=Rf(x)=\tilde{f}(p_0,x)$, $x\in X$.
 Notice that $\tilde{f}(\gamma p_0, \gamma x)=\tilde{f}(p_0,x)$. Since $p\mapsto \tilde{f}(p,x)$ is harmonic
 on $D$, its restriction to $V$ is $P$-harmonic
 by theorem \ref{fls}. Note  
 that \begin{equation*}\Phi_{p_0}(x)=\sum_{\xi\in \Gamma} \Phi_{p_0}(\xi x)\mu(\xi)\end{equation*} for all
 $x\in X$. I.e., $Rf$ belongs to $H(X, \Gamma, \mu)$.  
 In fact,
 \begin{align*}
 \Phi_{p_0}(x)-\sum_{\xi\in \Gamma} \Phi_{p_0}(\xi x)\mu(\xi)&=\tilde{f}(p_0,x)-
 \sum_{\xi\in \Gamma}\tilde{f}(p_0, \xi x)\mu(\xi)\\
 &=\tilde{f}(p_0,x)-
 \sum_{\xi\in \Gamma}\tilde{f}(\xi^{-1} p_0, x)\mu(\xi)\\
 &=0,
 \end{align*}
 by equation \ref{a} with $\gamma=e$, $\eta=\xi^{-1}$ and $\hat{f}(u)=\tilde{f}(u,x)$.
 This is because $\hat{f}$ is the restriction to $\Gamma\cdot p_0$ of the harmonic
 function $\tilde{f}(\cdot, x)$, using again theorem \ref{fls}.

 An equally straightforward manipulation gives the converse: 
 start with  $\Phi$ in $H(X, \Gamma, \mu)$ and define $\hat{g}:V\times X\rightarrow \mathbb{R}$ by
 $$\hat{g}(\gamma p_0, x):= \Phi(\gamma^{-1}x).$$ Then 
 the $P$-harmonic condition \ref{a}
 $$\hat{g}(\gamma p_0, x)=\sum_{\eta\in \Gamma} \hat{g}(\eta p_0, x)\mu(\eta^{-1}\gamma)$$
 holds. This is seen as follows:
 \begin{align*}
 \hat{g}(\gamma p_0, x)&=\Phi(\gamma^{-1}x)\\
 &=\sum_{\xi\in \Gamma}\Phi(\xi \gamma^{-1}x)\mu(\xi)\\
 &=\sum_{\eta\in \Gamma}\Phi(\eta^{-1}\gamma \gamma^{-1}x)\mu(\eta^{-1}\gamma)\\
 &=\sum_{\eta\in \Gamma} \hat{g}(\eta p_0, x)\mu(\eta^{-1}\gamma).
 \end{align*}

 By theorem \ref{fls}, $\hat{g}(\cdot, x)$ is the restriction to $V$ of a harmonic  function
 $ \tilde{g}(\cdot,x)$
 on $D$. The functions $p\mapsto \tilde{g}(p,\eta x)$ and $p\mapsto \tilde{g}(\eta^{-1}p, x)$ are
 both harmonic   and agree on $V$, so they must coincide on all of $D$,
 again by theorem \ref{fls}. Therefore $\tilde{g}$ 
 satisfies $\tilde{g}(\eta p, x )=\tilde{g}(p, \eta^{-1}x)$ for all $p, x, \eta$, hence there
 exists $g:M\rightarrow \mathbb{R}$ such that $g\circ \pi=\tilde{g}$.

 To conclude that $g\in H(M,\mathcal{F})$ we need to argue
  that $\tilde{g}$ is continuous on $D\times X$.
  Clearly $p\mapsto \tilde{g}(p,x)$ is continuous for each $x\in X$. 
  Also note that $x\mapsto \tilde{g}(\gamma p_0,x)$ is continuous for each 
  $\gamma\in \Gamma$. 
  We claim that
  $x\mapsto \tilde{g}(\cdot, x)$ is a continuous map from $X$ into the space
  of bounded harmonic functions on $D$.  By adding a positive constant 
  to $\tilde{g}$ if necessary,  we may
  assume without loss of generality that $\tilde{g}>0$. Now define
  $F(p,x)=\tilde{g}(p,x)/\tilde{g}(p_0,x)$. By the Harnack principle (see section \ref{probability})
  the space of non-negative harmonic functions $h$ on $D$ with the normalization
  $h(p_0)=1$ is compact in the topology of uniform convergence on compact subsets.
  Let $x_n$ be a sequence in $X$ converging to $x$. By passing to a subsequence we
  may assume that $F(\cdot, x_n)$ converges to a harmonic function $h$ on $D$.  We need
  to show that $h=F(\cdot, x)$.  But
   $h(\gamma p_0)=F(\gamma p_0, x)$ for all $\gamma\in \Gamma$  since $x\mapsto F(\gamma p_0,x)$
   is continuous.
  Therefore, $h(p)=F(p,x)$ for all $p\in D$ 
  by theorem \ref{fls}.  Multiplying $F$ back by the continuous function $\tilde{g}(p_0, \cdot)$
  implies that the sequence  of functions $\tilde{g}(\cdot, x_n)$ converges uniformly
  on compact sets to $\tilde{g}(\cdot, x)$. This proves that $\tilde{g}$ is continuous.

   Thus  we obtain $g\in H(M,\mathcal{F})$.
 It is a direct consequence of the definitions that $\Phi\mapsto g$ is the inverse operation
 to $R$. It
 is also clear that   leafwise constant functions on $M$ correspond to $\Gamma$-invariant functions on $X$ since 
 $$\tilde{g}(\gamma p_0, x)-\tilde{g}(p_0,x)=\Phi(\gamma^{-1}x)-\Phi(x) $$
 for all $\gamma\in \Gamma$ and $x\in X$.
\end{pf}

  Corollary \ref{dyn} now follows from theorem \ref{main}  and
 theorem \ref{discrete} applied to  $X=S^1$.  

\section{Discrete holomorphic functions}\label{holo}
It is interesting to note that the discretization  theorem allows one to define a notion of holomorphic
function in the discrete setting: let $\Gamma$ be a  group of covering transformations
of a simply connected,  transient K\"ahler manifold $D$
such that $D/\Gamma$ is  a compact manifold.
We call such $\Gamma$ a {\em transient K\"ahler group}.  Suppose  that $\Gamma$ acts
on a compact topological space $X$ by homeomorphisms. Given  a discretization measure $\mu$
on $\Gamma$ we say that a continuous $\Phi:X\rightarrow \mathbb{R}$ is $\mu$-holomorphic
if $R^{-1}\Phi$ is the real part of a  leafwise holomorphic function on the corresponding foliated $X$-bundle. (We may, of course, also consider complex-valued functions.) We recall that
$R$ is the restriction map defined immediately before theorem \ref{discrete}.

We illustrate the concept of $\mu$-holomorphic function by stating 
a discretized version of the following fact about foliations.

\begin{propo}[\cite{ghani1}]
Let $(M,\mathcal{F})$ a compact, connected foliated manifold with complex leaves. Suppose
that the closure of each leaf of $\mathcal{F}$ contains at most countably many minimal
sets. Then $\mathcal{F}$ has the holomorphic Liouville property.
\end{propo}

The next proposition follows immediately from the previous 
one   and  theorem \ref{discrete}.

\begin{propo}
Let  a transient K\"ahler group $\Gamma$ act by homeomorphisms on a compact
topological space $X$, and let $\mu$ be a discretization measure on $\Gamma$.
Suppose that the closure of each $\Gamma$-orbit contains at most
a countable number of minimal sets. Then every continuous $\mu$-holomorphic function
on $X$ is $\Gamma$-invariant.
\end{propo}

We give now  an example of a $\Gamma$-action 
that admits non-trivial $\mu$-holomorphic functions. 
 The example is a modified version of
the one  shown in  \cite{ghani1} immediately after theorem 1.16.
To make the construction more transparent, 
we let the space $X$   be a manifold with boundary, but we can 
also obtain an action on  a manifold without boundary by doubling. 
Representing
an element of $\mathbb{R}P^3$ as $[u, v]$, $u, v\in \mathbb{C}$ not both
$0$, let $X$ be the subset  of all $[u, v]$ such that $|u|\geq |v|$. Thus $X$ is a solid torus
(it is doubly covered by $\{(u,v)\in \mathbb{C}^2: |u|=1, |v|\leq 1\}$.)
The boundary of $X$  consists of all $[e^{i\theta}, e^{i\varphi}]$, $\theta, \varphi\in \mathbb{R}$,
so $\partial X$ is homeomorphic to a $2$-torus. Notice that   $U(1)=\{e^{i \xi}: \xi\in \mathbb{R}\}$ acts
on $X$ by $\omega [u,v]=[ \omega u,  \omega v]$ leaving  the boundary invariant and
having circle orbits. This defines a Seifert fibration in $X$.
The $U(1)$-action foliates the boundary by circles and the space of leaves of
$\partial X/ U(1)$ is also a circle.

It is well known that the group of isometries of $\mathbb{D}$ is isomorphic to
$PSL(2, \mathbb{R})$.  It is somewhat more convenient  to
 use the isomorphic representation of it as $G=PSU(1,1)$,
the group of all $2\times 2$ complex valued matrices of the form
$\gamma=\left(\begin{array}{cc}{\alpha} & \beta \\ \overline{\beta} & \overline{\alpha}\end{array}\right)$
modulo the center, $\pm I$, where $|\alpha|^2 - |\beta|^2=\det \gamma=1$. 
The action on $\mathbb{R}P^3$  defined by 
$\gamma [u, v]=[\alpha u + \beta \overline{v}, \alpha v+ \beta \overline{u}]$
is easily shown to have the following properties:
\begin{enumerate}
\item  $X$  and $\partial X$ are invariant sets. In fact, writing  $r=\alpha u + \beta \overline{v}$
and $s=\alpha v+ \beta \overline{u}$ so that $\gamma [u,v]=[r,s]$, 
 then
it is easily checked that $$|r|^2 - |s|^2=(|\alpha|^2- |\beta|^2)(|u|^2 - |v|^2)=|u|^2-|v|^2\geq 0;$$
\item For each $\gamma\in G$ and $[u,v]\in \partial X$, one has
$\gamma [u,v]=[\omega u, \omega v]$
for some $\omega\in  U(1)$.
(Observe that $(\alpha u + \beta \overline{v})/u=(\alpha v + \beta \overline{u})/v=\alpha + \beta \overline{uv}$ if $|u|=|v|=1$. Therefore, $\omega=(\alpha + \beta \overline{uv})/|\alpha + \beta \overline{uv}|$.)
 In particular, the $U(1)$-orbits in $\partial X$
are also invariant; 
\item For any $[u,v]\in X$, 
   $\gamma [u,v]$ approaches
the torus boundary as $\gamma\rightarrow \infty$ in $G$. 
In fact, by the formula and notation of item (1) we see that 
$$1- \frac{|s|^2}{| r|^2}=\frac{|u|^2-|v|^2}{|r|^2}.$$
As $\gamma\rightarrow \infty$ in $G$,  it is easily checked that $|r|\rightarrow \infty$, 
hence the claim.
Therefore, any minimal set
for the action of any non-compact subgroup of $G$ on $X$ is contained in one
$U(1)$-orbit in $\partial X$.  
\end{enumerate}

Now let $(M,\mathcal{F})$ be the foliated $X$-bundle 
over a compact $\mathbb{D}/\Gamma$
associated to the given action restricted to $\Gamma$. 
Then $M$ is a $5$-manifold and   $\mathcal{F}$ has (real) codimension $3$.
 Define $f\in H(M,\mathcal{F})$ such that $\tilde{f}:\mathbb{D}\times X\rightarrow
 \mathbb{C}$ is given by
 $$\tilde{f}(z,[u,v]):=({\overline{u}z - v})/(u - {\overline{v}z  }). $$
 This definition only  makes sense {\em a priori}  in the interior of $X$, but as 
 $[u, v]$ approaches a boundary point $[e^{i\theta}, e^{i\varphi}]$ the function
 $\tilde{f}(\cdot, [u,v])$ converges to the constant   $-e^{i(\varphi-\theta)}$
 uniformly on compact subsets of $\mathbb{D}$.
 Notice how  this limit   is the same along the $U(1)$-orbits in $\partial X$.
A straightforward calculation shows that $\tilde{f}$ is $\Gamma$-invariant and
so defines a continuous function $f$ on $M=(\mathbb{D}\times X)/\Gamma$ which
   is  leafwise holomorphic and
non-constant on all interior leaves. By the discretization theorem we obtain 
a continuous  $\mu$-harmonic function on $X$ which is not constant on 
interior $\Gamma$-orbits. All orbits
 accumulate
on the boundary of $X$.

\begin{small}

DEPARTMENT OF MATHEMATICS, FLORIDA STATE UNIVERSITY,
TALLAHASSEE, FL 32306 USA

\textit{E-mail address}: fenley@math.fsu.edu
\bigskip

DEPARTMENT OF MATHEMATICS,
WASHINGTON UNIVERSITY, ST.~LOUIS, MO 63130 USA

\textit{E-mail address}: feres@math.wustl.edu
\bigskip

DEPARTMENT OF MATHEMATICS, EASTERN ILLINOIS UNIVERSITY, 
CHARLESTON, IL 61920 USA

\textit{E-mail address}: forty2@math.northwestern.edu 

 \end{small}
 
\end{document}